\documentclass[a4paper,11pt,reqno]{amsart}
\usepackage{amsthm}
\usepackage{amssymb}
\usepackage{latexsym}
\usepackage{multicol}
\usepackage{verbatim,enumerate}
\usepackage{accents}
\usepackage[usenames]{color}
\usepackage{hyperref}
\usepackage{amsmath, amscd}
\usepackage{soul}
\usepackage{tikz}
\usepackage{youngtab}
\usepackage{epic, eepic}
\usepackage{float}

\usetikzlibrary{matrix,arrows,decorations.pathmorphing}
\advance\textwidth by 1.2in \advance\oddsidemargin by -.6in
\advance\evensidemargin by -.6in

\numberwithin{equation}{section} 
\theoremstyle{plain}
\newtheorem*{prop}{Proposition}
\newtheorem{thm}{Theorem}
\newtheorem*{lem}{Lemma}
\newtheorem*{cor}{Corollary}

\theoremstyle{definition}
\newtheorem*{rem}{Remark}
\newtheorem*{example}{Example}
\newtheorem*{defn}{Definition}

\theoremstyle{remark}

\newcommand{\lie}[1]{\mathfrak{#1}}

\newcommand\bc{\mathbb C}
\newcommand\bn{\mathbb N}
\newcommand\bz{\mathbb Z}

\newcommand{\ev}{\operatorname{ev}}
 \newcommand\bob{\bold b}            
          
\newcommand{\qbinom}[2]{\genfrac[]{0pt}0{#1}{#2}}

\begin{document}

\makeatletter
\def\section{\def\@secnumfont{\mdseries}\@startsection{section}{1}%
  \z@{.7\linespacing\@plus\linespacing}{.5\linespacing}%
  {\normalfont\scshape\centering}}
\def\subsection{\def\@secnumfont{\bfseries}\@startsection{subsection}{2}%
  {\parindent}{.5\linespacing\@plus.7\linespacing}{-.5em}%
  {\normalfont\bfseries}}
\makeatother

\title[A combinatorial formula for graded multiplicities in excellent filtrations]{A combinatorial formula for graded multiplicities in excellent filtrations}
\author{Rekha Biswal}
\address{Universit\'e Laval, D\'e{}partment de math\'e{}matiques et de Statistique, Qu\'e{}bec, QC, Canada}
\email{rekha.biswal.1@ulaval.ca}
\thanks{}

\author{Deniz Kus}
\address{University of Bochum, Faculty of Mathematics, Universit{\"a}tsstr. 150, 44801 Bochum, 
Germany}
\email{deniz.kus@rub.de}
\thanks{}
\subjclass[2010]{}
\begin{abstract}
A filtration of a representation whose successive quotients are isomorphic to Demazure modules is called an excellent filtration. In this paper we study graded multiplicities in excellent filtrations of fusion products for the current algebra $\mathfrak{sl}_2[t]$. We give a combinatorial formula for the polynomials encoding these multiplicities in terms of two dimensional lattice paths. Corollaries to our main theorem include a combinatorial interpretation of various objects such as the coeffficients of Ramanujan's fifth order mock theta functions $\phi_0, \phi_1, \psi_0,  \psi_1$, Kostka polynomials for hook partitions and quotients of Chebyshev polynomials. We also get a combinatorial interpretation of the graded multiplicities in a level one flag of a local Weyl module associated to the simple Lie algebras of type $B_n \text{ and } G_2$. 
\end{abstract}

\maketitle

\section{Introduction}
In this paper, we are interested in objects in the category of finite-dimensional $\bz$--graded modules for the current algebra $\mathfrak{sl}_2[t]$ which admit a filtration 
whose successive quotients are isomorphic to stable Demazure modules. These types of filtrations are called excellent filtrations \cite{Jo03,J06} but are also known in the literature as Demazure flags \cite{BCSV15,CSSW14}. Demazure modules occur in highest weight integrable modules of the affine Lie algebra $\widehat{\mathfrak{sl}_2}$ and are called stable if they admit an action of the maximal parabolic subalgebra $\mathfrak{sl}_2[t]$; the integer given by the action of the canonical central element is called its level. The great interest in Demazure modules is due to several reasons. They appear as classical limits of a family of irreducible representations of the quantum affine algebra \cite{BCMo15,CP01}, their graded characters are specializations of Macdonald polynomials \cite{I03} and satisfy certain functional relations which are known as $Q$--systems \cite{CV13, KV14}, just to name a few.\par
An important result due to Naoi states that if $m \geq m' \geq 1$, then a stable Demazure module of level $m'$ admits a filtration such that the successive quotients are isomorphic to level $m$--Demazure modules. In fact Naoi proves this result for an affine Lie algebra associated to a simply laced simple Lie algebra. Later, this result was greatly extended for a wider class of modules for $\mathfrak{sl}_2[t]$, which are indexed by partitions and are known in the literature as fusion products \cite{FL99}.
 They generalize many known families of representations for $\mathfrak{sl}_2[t]$, e.g. stable Demazure modules \cite[Theorem 2]{CV13}, local Weyl modules \cite{CP01} and tuncated Weyl modules \cite[Theorem 4.3]{KL15}. Under a suitable condition on the partition, fusion products admit a level $m$--Demazure flag \cite[Theorem 3.3]{CSSW14} and recurrence relations were established in the same paper. A closed form solution of these recurrences was however, only obtained in very special cases (see Remark~\ref{rem21}). It turned out that in these special cases there is a beautiful link between these modules and number theory and combinatorics: numerical multiplicities are closely related to Chebyshev polynomials \cite[Corollary 1.3]{BCSV15}, several specializations of the generting series associated to the graded multiplicities specialize to Ramanujan's fifth order mock theta functions \cite[Theorem 1.6.]{BCSV15}, certain weighted versions of the generating series give Carlitz q--Fibonacci polynomials \cite[Proposition 2.5.3]{BCK16} and are limits of hypergeometric series \cite[Section 2.5.4]{BCK16}. One of the motivations of the present paper is to complete the results of \cite{BCSV15, CSSW14} by finding a formula for the graded multiplicities in a level $m$--Demazure flag of an arbitrary fusion product and to further understand the deep and unexpected link between the theory of Demazure flags and combinatorics. The main ingredients in our study are two dimensional lattice paths.\par
A lattice path is a path in a lattice in some euclidean vector space. They have a long history and have entered many fields of mathematic, computer science and physics; for a survey of results in the enumeration of lattice paths we refer the reader to \cite{Kra15}. A diagonal lattice path which consist only of up--steps $(1,1)$ and down--steps $(1,-1)$ and which do not pass below the $x$--axis is called a Dyck path. There are several statistics on the set of Dyck paths and we are interested in the comajor statistics  (see Section~\ref{section35}) which was studied first by MacMahon \cite{MacM60} in his interpretation of the $q$--Catalan numbers. Our strategy for determining graded multiplicities in excellent filtrations of fusion products is as follows.\par Let $\xi$ a partition and $V(\xi)$ the fusion product associated to $\xi$. We fix a level $m$--Demazure module $\tau_p^* D(m,n)$ of heighest weight $n$ whose degree is $p$ (see Section~\ref{section2} for the precise definitions). Further, set
$$\mathcal{V}_{n}^{\xi\rightarrow m}(q)=\sum_{p\ge 0}\ [V(\xi):\tau_p^* D(m,n)] \hspace{0,04cm} q^p,$$
where  $[V(\xi):\tau_p^* D(m,n)]$ is the multiplicity of $\tau_p^* D(m,n)$ in a level $m$--Demazure flag of $V(\xi)$. In the first step, we reduce the computation of graded multiplicities in excellent filtrations of fusion products to the computation of graded multiplicities in local Weyl modules. To be more precise, we show that there exists a polynomial in finitely many variables such that a suitable evaluation of that polynomial at elements of the form $\mathcal{V}_{\bullet}(q)^{(1,1,\dots,1)\rightarrow \bullet}$ gives the desired polynomial $\mathcal{V}_n(q)^{\xi\rightarrow m}$ (see Poposition~\ref{reductionpr}). In the second step we give a combinatorial formula for $\mathcal{V}_{n}^{\xi\rightarrow m}(q)$ in terms of Dyck paths when $\xi$ is a hook partition. Let $m,m'\in\bn$, $s\in\bz_+$ such that $m\geq m'$ and $\xi=(m',1^s)$ a hook partition. Then we have (see Theorem~\ref{mainthm1}) 
$$\mathcal{V}^{\xi\rightarrow m}_{n}(q)=\sum_{P} q^{\text{comaj(P)}},$$
where the sum runs over the set of admissible Dyck paths (see Definition\ref{maindef}) from the origin $(0,0)$ to $(s+m',n)$ which start with $m'$ up--steps and never cross the line $y=\max\{m-1,n\}$.\par
The appearance of Dyck paths in this set up is quite unexpected and intriguing. A number of consequences flow from our new combinatorial formula. We briefly summarize the consequences below:
\begin{itemize}
\item[(i)] We find new combinatorial interpretations of Ramanujan's fifth order mock theta functions $\phi_1, \phi_0, \psi_1, \psi_0$ in terms of Dyck paths and their comajor index (see Corollary~\ref{mockcoef}). 
\item[(ii)] We find a combinatorial interpretation of the quotients of powers of two Chebyshev polynomials in terms of Dyck paths of bounded height generalizing the results of Gessel and Xin \cite{GX05}  (see Proposition~\ref{numerical}). 
\item [(iii)]The powers of $q$ which appear with non--zero coefficients in the polynomial encoding the graded multiplicities in Demazure flags of local Weyl modules form an interval of consecutive integers (see Corollary~\ref{conint}). 
\item [(iv)]The graded multiplicities of irreducible modules in local Weyl modules are principal specializations of Schur functions (see Corollary~\ref{princp}). 
\item [(v)] We give formulae for the graded multiplicities of an irreducible module in a Demazure flag of a local Weyl module (see Corollary~\ref{KrM}).
\end{itemize}
Our paper is organized as follows. In Section 2 we recall the notion of excellent filtrations. In Section 3 we define our combinatorial model which is the crucial object in our main theorem. In Section 4 we reduce the computation of graded multiplicities in excellent filtrations of fusion products to the computation of graded multiplicities in local Weyl modules and state our main theorem. We also discuss several corollaries. In Section 5 we give a proof of our main theorem by using a different recursion which is obtained from a representation theoretical result proved in Section 6.  

\textit{Acknowledgements: R.B. thanks Francois Bergeron, Ira Gessel and Dennis Stanton for many fruitful discussions about Chebyshev polynomials which was the starting point of this work. D.K. thanks Christian Krattenthaler for many helpful discussions on the combinatorics of Dyck paths and for pointing him to formula \eqref{b}. R.B. also gratefully acknowledges the funding received from  NSERC discovery grant of her Postdoc supervisor Michael Lau and  Universit\'e Laval for the hospitality.}

\section{Excellent filtrations and graded multiplicities}\label{section2}
In this section we recall the notion of excellent filtrations (Demazure flags) \cite{Jo03,J06} and set up the notation needed in the rest of the paper.

\subsection{}
Throughout this paper we denote by $\mathbb{C}$ the field of complex numbers and by $\mathbb{Z}$ (resp. $\mathbb{Z}_{+}$, $\mathbb{N}$) the subset of integers (resp. non-negative, positive integers). Given $n,m\in\bz$, set
\begin{gather*}\qbinom{n}{m}_q=\frac{(1-q^n)\cdots(1-q^{n-m+1})}{(1-q)\cdots (1-q^m)}, \ \ \  n\geq m>0 ,\\  \\ \qbinom{n}{0}_q=1,\ \  n\geq 0,\ \ \qbinom{n}{m}_q= 0,\ \  \text{otherwise}.\end{gather*}
Further, we introduce the $q$--Pochammer symbol $(a;q)_{n}=\prod_{i=0}^{n-1}(1-aq^{i}).$
\subsection{}  Let $\lie{sl}_2[t]\cong\lie{sl}_2\otimes \bc[t]$ the Lie algebra of two by two matrices of trace zero with entries in the algebra $\bc[t]$ of polynomials with complex coefficients. The degree grading of $\bc[t]$ defines a natural grading on $\lie{sl}_2[t]$.  A  finite-dimensional $\bz$-graded $\lie {sl}_2[t]$--module is a $\bz$--graded vector space admitting a compatible graded action of $\lie{sl}_2[t]$:
  $$ V=\bigoplus_{k\in\bz} V[k],\qquad (a\otimes t^r)V[k]\subset V[k+r]\ \  a\in \lie{sl}_2,\ \ r\in\bz_+.$$  Given a $\bz$--graded space $V$ let $\tau_p^* V$ the graded vector space whose $r$--th graded piece is $V[r+p]$.
\subsection{}\label{section23} The category of finite-dimensional $\bz$--graded $\lie {sl}_2[t]$--modules was the central subject of many recent papers (see, for example \cite{BC15,BCM12,BCSV15,CSSW14,CV13,FL99,KL15}). There is a well--known family of objects in that category which is of particular interest, namely the subclass of fusion products. We recall their description in terms of generators and relations from \cite[Section 6]{CV13}; for a more traditional definition we refer the reader to \cite{FL99} (see also Section~\ref{551}). Let $x,h,y$ be the standard basis of $\lie{sl}_2$ and set $u^{(r)}:=\frac{1}{r!} u^r$ for $u\in \lie{sl}_2[t],\hspace{0,03cm} r\in\bz_+$. For a partition $\xi=(\xi_1\geq \xi_2\geq \cdots \geq \xi_\ell)$ we set 
$$|\xi|_i:=\sum_{j=i}^{\ell}\xi_j,\ 1\leq i \leq \ell, \ \ |\xi|:=|\xi|_1.$$
The fusion product associated to $\xi$ is the $\lie {sl}_2[t]$--module $V(\xi)$ generated by an element $v_{\xi}$ with defining relations:
\begin{gather}\label{demreldes}(x\otimes \bc[t])v_{\xi}=0,\ \ (h\otimes f) v_{\xi}= |\xi| f(0)v_{\xi},\ \  \  (y\otimes 1)^{|\xi|+1}v_{\xi}=0,\\ \label{demrel2a}  
(x\otimes t)^{(p)}(y\otimes 1)^{(r+p)}v_{\xi}=0,\ \ r,p\in \bn,\ r+p\geq 1+rk+\sum_{j\geq k+1}\xi_j\ \ \text{for some $k\in \bn$.}
\end{gather}  
It turns out that many other well--known families of representations belong to the class of fusion products. For example, Demazure modules occur in irreducible integrable representations of the affine Lie algebra $\widehat{\mathfrak{sl}_2}$ and are parametrized by tuples $(m,n)\in \bn\times \bz_+$, where the integer $m$ is called the level. We denote such a module by $D(m,n)$. If $n_0,n_1\in\bz_+$ are such that $n_0<m$ and $n=n_1m+n_0$, then the fusion product $V(\xi(m,n))$ associated to the partition $\xi(m,n):=(m^{n_1},n_0)$ is isomorphic to $D(m,n)$ (see \cite[Theorem 2]{CV13}). Hence Demazure modules can be categorized into the family of fusion products, but the class of fusion products is generically much bigger. Nevertheless, there is a beautiful result saying that fusion products admit a filtration by  Demazure modules under a suitable condition on the partition. The following proposition was proved in \cite[Theorem 3.3]{CSSW14}.
\begin{prop}\label{existflag}
Let $\xi=(\xi_1\geq \xi_2\geq \cdots \geq \xi_{\ell})$ a partition and $m\in\bn$. The module $V(\xi)$ admits a filtration of level $m$,  i.e., there exists a decreasing sequence of graded submodules
$$0=V_0\subset V_1\subset \cdots V_{k-1}\subset V_k=V(\xi)$$
 such that $$V_i/V_{i-1}\cong \tau_{p_i}^*D(m, n_i),\ \ (p_i,n_i)\in \bz_+\times \bz_+,\ \ 1\leq i\leq k$$ 
 if and only if $m\geq \xi_1$.
\hfill\qed
\end{prop}
These types of filtrations are called excellent filtrations \cite{Jo03,J06} but are also known in the literature as level $m$--Demazure flags \cite{BCSV15,CSSW14}. The aim of the present paper is to give a combinatorial formula for the graded multiplicities in excellent filtrations.
\begin{rem}\label{rem23} We discuss two further specializations of the partition $\xi$. 
\begin{enumerate}
\item The specialization $\xi_i=1$, $1\leq i\leq \ell$ leads to a representation which is isomorphic to the local Weyl module $W_{\text{loc}}(|\xi|)$ which is the module generated by an element $w_{|\xi|}$ subject to the relations \eqref{demreldes}. The interest in local Weyl modules is its connections with quantum affine algebras \cite{CP01}, the theory of Macdonald polynomials \cite{I03}, $q$--Whittaker functions \cite{BF14a}, $Q$--systems \cite{CV13, KV14} and more recently with hypergeometric series \cite{BCK16}.
\item Let $N\in \bn$. The truncated Weyl module $W_{\text{loc}}(N,n)$ is a quotient of $W_{\text{loc}}(n)$ by the additional relation $(\mathfrak{sl}_2\otimes t^N\bc[t])w_{n}=0$. The special choice $\xi=((d+1)^j,d^{N-j})$ yields an isomorphism $V(\xi)\cong W_{\text{loc}}(N,n)$, where $d\in \bz_+$ and $j<N$ are such that $n=dN+j$ (see \cite[Theorem 4.3]{KL15}).

\end{enumerate}
\end{rem}
We will secretly assume in the rest of the paper that $m\geq \xi_1$ whenever we talk about level $m$--Demazure flags of $V(\xi)$. 
\subsection{}The number of times a particular level $m$--Demazure module appears as a quotient in a level $m$--flag is independent of the choice of the flag. We encode these multiplicities in a polynomial
\begin{equation}\label{poly}
\mathcal{V}_{n}^{\xi\rightarrow m}(q):=\sum_{p\ge 0}\ [V(\xi):\tau_p^* D(m,n)] \hspace{0,04cm} q^p,\end{equation}
where  $[V(\xi):\tau_p^* D(m,n)]=\sharp\{1\leq i\leq k : V_i/V_{i-1}\cong \tau_p^* D(m,n)\}$. It is known that 
\begin{equation}\label{initial12}
\mathcal{V}_{s}^{\xi(m',s)\rightarrow m}(q)=1,\quad \mathcal{V}_{n}^{\xi(m,s)\rightarrow m}(q)=\delta_{s,n},\ \ \mathcal{V}_{n}^{\xi\rightarrow m}(q)=0,\ \ \text{if }\  |\xi|-n\notin2\bz_+.\end{equation} Moreover, for $m\geq m'\geq \xi_1$ we have 
\begin{equation}\label{mat}\mathcal{V}_{n}^{\xi\rightarrow m}(q)=\sum_{p\geq 0} \mathcal{V}_{p}^{\xi\rightarrow m'}(q)\ \mathcal{V}_{n}^{\xi(m',p)\rightarrow m}(q).\end{equation}
For convinience, set
\begin{equation}\label{initial1}
\mathcal{V}_{n}^{\xi\rightarrow m}(q)=0,\  \ {\rm{if}} \ \ |\xi|<0\ \   {\rm{or}}\  n<0.
\end{equation}
The following lemma will be needed later; for a proof see  \cite[Lemma 3.8]{CSSW14}.
\begin{lem}\label{ausc}Let $\xi=(\xi_1\geq \xi_2\geq \cdots \geq \xi_{\ell})$ a partition and set $\xi'=(\xi_2\geq \xi_3\geq \cdots \geq \xi_{\ell})$. We have
$$\mathcal{V}_{n}^{\xi\rightarrow \xi_1}(q)=q^{(|\xi|-n)/2}\hspace{0,04cm} \mathcal{V}_{n-\xi_1}^{\xi'\rightarrow \xi_1}(q).$$
\hfill\qed
\end{lem}
We will freuquently deal with the polynomials \eqref{poly} when $\xi=\xi(m',s)$, i.e. $V(\xi)$ itself is a level Demazure module of level $m'$. So we abbreviate
$$\mathcal{V}_{s,n}^{m'\rightarrow m}(q):=\mathcal{V}_{n}^{\xi(m',s)\rightarrow m}(q)$$
and define the associated generating series by
$$\mathcal{A}^{m'\rightarrow m}_n(x,q)=\sum_{s\geq 0}\mathcal{V}_{s,n}^{m'\rightarrow m}(q)\,x^{(s-n)/2},\ \ \ n\ge 0.$$ 
\begin{rem}\label{rem21}
The proof of Proposition \ref{existflag} uses the short exact sequence of \cite[Theorem 5]{CV13}. It has the advantage that one can derive recursive formulas for $\mathcal{V}_{s,n}^{m'\rightarrow m}(q)$ (see for example \cite[Theorem 2.2]{CSVW14}) but closed formulas were established only in the case of $m'=1,m=2$ \cite[Theorem 3.3]{CSSW14} and $m'=1,m=3$ \cite[Section 1.6]{CSVW14}. Our approach uses a different short exact sequence (see Section \ref{section6}) which enables us to give a combinatorial formula in the most general setting.
\end{rem}
\section{Combinatorics of Dyck paths}\label{section3}
In this section we introduce our combinatorial model and certain combinatorial statistics which will be needed to describe the polynomials \eqref{poly}.
\subsection{} \label{section31} A Dyck path is a diagonal lattice path from the origin $(0,0)$ to $(s,n)$ for some non--negative integrs $s,n\in \bz_+$,
such that the path never goes below the x--axis. We encode such a path by a 01--word, where $1$ encodes the up--steps and $0$ the down--steps. For a path $P$ we set $$\text{supp}(P):=\{z\in \bz_+^2 : z \text{ is a point on $P$}\}.$$
The length of a point $z$ on a path $P$ is defined to be its $x$--coordinate and its $y$--coordinate is called the height; we will frequently use the notation $z=(\ell_z(P),\text{ht}_z(P))$. We denote by $\mathcal{D}_{n}$ the set of Dyck paths that end at height $n$ and by $\mathcal{D}_n(s)$ the subset of paths $P$ with $\ell(P)=s$, where $\ell(P)$ denotes the length of the endpoint of $P$. Obviously (compare with \eqref{initial12}),
$$\mathcal{D}_{n}=\bigcup_{s\geq 0} \mathcal{D}_n(s),\ \ |\mathcal{D}_s(s)|=1, \ \ |\mathcal{D}_n(0)|=\delta_{n,0}, \ \ \mathcal{D}_n(s)=\emptyset, \ \text{ if } s-n\notin 2\bz_+.$$
Further, let $\mathcal{D}^{m}_n$ the subset of paths which do not cross the line $y=m$, i.e.
$$\mathcal{D}^{m}_n=\{P\in \mathcal{D}_{n}:  \text{ht}_z(P)\leq m \text{ for all } z\in \text{supp}(P)\},\ \ \mathcal{D}_n^{m}(s):=\mathcal{D}_n^m\cap \mathcal{D}_n(s).$$ 
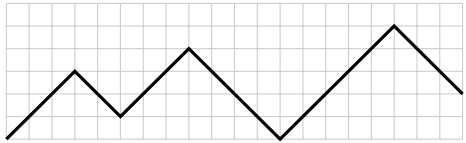
\begin{figure}[H]
\begin{tikzpicture}[scale=0.3]
\draw [step=1,thin,gray!40] (0,0) grid (20,6);
\draw [very thick] (0,0) -- (3,3) -- (5,1) -- (8,4)--(12,0)--(17,5)--(20,2);
\end{tikzpicture}
\caption{A Dyck path $P\in \mathcal{D}_2^5(20)$}
\label{fig1}
\end{figure}
In analogy to \eqref{initial1} we set
\begin{equation}\label{initial2}
\mathcal{D}_n(s)=0,\  \ {\rm{if}} \ \ s<0\ \   {\rm{or}}\  n<0.
\end{equation}
\begin{rem}
In the literature a Dyck path mostly ends at height $0$ and a diagonal lattice path that ends at an arbitrary height is often called a generalized Dyck path. It is well--known that the cardinality of $\mathcal{D}_{0}(2k)$ is given by the $k$--th Catalan number; for further combinatorial interpretations of the Catalan numbers we refer to the book of Stanley \cite[pg. 219-229]{S99}.
\end{rem}
\subsection{}The generating function $\mathcal{G}_n^m(x)=\sum_{s\geq 0}|\mathcal{D}^m_{n}(s)| x^{(s-n)/2}$ has been intensively studied in \cite{GX05}. It turns out that $\mathcal{G}_n^m(x)$ is a rational function which can be expressed in terms of Chebyshev polynomials. We shall explain this connection in more detail. The Chebyshev polynomials of the second kind are defined by the recurrence relation 
$$U_n(x)=2x U_{n-1}(x)-U_{n-2}(x),\ \ n\geq 2,$$
with initial data $U_0(x)=1$ and $U_1(x)=2x$. Define polynomials
$$p_n(x):=x^{n/2}U_n((2\sqrt{x})^{-1})=\sum_{s=0}^{\lfloor\frac{k}{2}\rfloor} (-1)^s \qbinom{k-s}{s}.$$
The first equality of the next lemma follows from \cite[Lemma 4.2]{GX05} and the second equality has been proved in \cite[Corollary 1.3]{BCSV15}.
\begin{lem}\label{Ges} We have $\mathcal{G}^m_n(x)=0$ if $m<n$ and otherwise
$$\mathcal{G}^m_n(x)=\frac{p_{m-n}(x)}{p_{m+1}(x)}=\mathcal{A}_{n}^{1\rightarrow m+1}(x,1).$$ In particular,
$$\mathcal{V}_{s,n}^{1\rightarrow m+1}(1)=|\mathcal{D}^m_n(s)|,\ \ \forall n\leq m.$$
\hfill\qed
\end{lem}
This result suggests a deeper connection between graded multiplicities in Demazure flags and certain Dyck path statistics and the motivation of the present paper is to determine this connection. As a byproduct we will generalize the above lemma and find a combinatorial model (in terms of Dyck paths) whose generating series is $\mathcal{A}^{1\rightarrow m}_n(x,1)$ (even for $n\geq m$) (see Proposition \ref{numerical}).
\subsection{}\label{admissibledef} A point $(z_1,z_2)\in \text{supp}(P)$ is called a peak (resp. valley) of the path if $(z_1\pm 1,z_2-1)\in \text{supp}(P)$ (resp. $(z_1\pm 1,z_2+1)\in \text{supp}(P)$). By convention, we will also call the endpoint of a path which is immediately preceded by a down--step a valley.  For example, the path in Figure \ref{fig1} has three peaks and three valleys.
Given a pair of non-negative integers $(a,b)\in \bz^2_+$, we say that $P\in \mathcal{D}_n$ is $(a,b)$--admissible if and only if $P$ satisfies the following property
\begin{equation}\label{admissible} \text{$P$ has a peak $z$ at height $b$} \Rightarrow \text{$\text{ht}_{z'}(P)>a$ for all $z'\in \text{supp}(P)$ with $\ell_{z'}(P)\geq \ell_{z}(P)$}.\end{equation}
With other words, once the path has a peak at height $b$, the remaining subsequent path is strictly above the line $y=a$ (see Figure~\ref{fig55} for a pictorial illustration).
\begin{figure}[H]
\begin{tikzpicture}[scale=0.3]
\draw [step=1,thin,gray!40] (0,0) grid (16,6);
\draw [dotted] (0,0) -- (3,3);
\draw [very thick] (3,3) -- (4,4) -- (5,3);
\draw [dotted, red,very thick] (4,1) -- (16,1);
\node [left=1pt] at (0,1) {\scriptsize $y=a$};
\node [left=1pt] at (0,4) {\scriptsize $y=b$};
\end{tikzpicture}
\caption{}
\label{fig55}
\end{figure}
\subsection{}For $n,m\in \bz_+$ let $n_0,n_1\in\bz_+$ such that $n_0<m$ and $n=m n_1+n_0$. \textit{In what follows we set $N=\max\{m-1,n\}$.} We are interested in a suitable subset of $\mathcal{D}^N_n$ which we will call admissible Dyck path. If $n<m$ we set $A(m,n)=\emptyset$ and otherwise define
$$A(m,n):=\{(i_1,m),(i_2,m+1),\dots,(i_{n-m+1},n)\}\subset \bz_+^2$$
where $i_1<\cdots<i_{n-m+1}$ is the natural ordering of the set 
$$\{0,\dots,n\}\backslash \{pn_1+n_0+\min\{0,(p-1)-n_0\},\ 1\leq p\leq m\}.$$
Equivalently, 
\begin{align}A(m,n)=\notag &\{(p(n_1+1)+r,m+pn_1+r),\ 0\leq r<n_1,\ 0\leq p\leq n_0\}\ \cup &\\& \notag \{(pn_1+(n_0+1)+r,m+p(n_1-1)+(n_0+1)+r),\ 0\leq r<n_1-1,\ n_0< p<m\}.\end{align}
The following is straightforward:
\begin{equation}\label{1}
(a,b)\in A(m,n)\backslash\{(0,m)\} \Rightarrow \exists \ \tilde a<a: \ \  (\tilde a,b-1)\in A(m,n).
\end{equation}
\begin{defn}\label{maindef}
We call a path $P\in \mathcal{D}_n^N$ admissible iff $P$ is $(a,b)$--admissible for all $(a,b)\in A(m,n)$. We denote by $\mathcal{D}_{m,n}$ the set of admissible Dyck paths and set $\mathcal{D}_{m,n}(s):=\mathcal{D}_n^N(s)\cap \mathcal{D}_{m,n}$.
\end{defn}
 Both sets will play a major role in the description of the generating series and graded multiplicities respectively, see Theorem \ref{mainthm1}. 
\begin{example}We have $A(n,n)=\{(0,n)\}$ and hence $\mathcal{D}_{n,n}$ consists of all Dyck paths $P\in \mathcal{D}^n_{n}$ which do not return to the $x$--axis after a peak of height $n$. The green path in Figure \ref{fig3} is admissible and the red path violates the condition and is not admissible.
\begin{figure}[H]
\begin{tikzpicture}[scale=0.3]
\hspace{-3cm}
\draw [step=1,thin,gray!40] (0,0) grid (16,5);
\draw [green,very thick] (0,0) -- (4,4) -- (6,2) -- (7,3) -- (8,2) -- (10,4) -- (13,1) -- (16,4);
\end{tikzpicture}
\begin{tikzpicture}[scale=0.3]
\draw [step=1,thin,gray!40] (0,0) grid (16,5);
\draw [red, very thick] (0,0) -- (4,4) -- (6,2) -- (7,3) -- (10,0) -- (12,2) -- (13,1) -- (16,4);
\end{tikzpicture}
\caption{$n=4$}
\label{fig3}
\end{figure}
\end{example}
\subsection{}\label{section35} The major statistics of a Dyck path was studied first by MacMahon \cite{MacM60} in his interpretation of the $q$--Catalan numbers. Let $P=a_1\cdots a_{s}$, $a_i\in\{0,1\}$ a Dyck path of length $s$. The major and comajor index are defined by 
$$\text{maj}(P)=\sum_{\substack{1\leq i < s,\\ a_i>a_{i+1}}} i,\ \ \ \text{comaj}(P)=\sum_{\substack{1\leq i < s,\\ a_i>a_{i+1}}} (s-i).$$
\begin{rem}\label{rem1}
The comajor index can also be defined for standard Young tableaux and play an important role in the expression of the character of the space of harmonics (the vector space spanned by the Vandermonde determinant and its partial derivatives of all orders; see the work of Haiman \cite{Hai94}). The definition is essentially the same in the sense that there is a bijection from $\mathcal{D}_n(n+2k)$ to the set of standard Young tableaux of shape $\lambda=(n+k,k)$, which preserves the comajor statistics (see Section~\ref{section33}). However, the tranlation of the admissibility conditions seem to be quite technical and hence we prefer to work with the notion of Dyck paths instead of two row partitions. 
\end{rem}

\section{The main results}\label{section4}
In this section we summarize the main results of the paper. We have arranged it so that this section can be read essentially independently of the representation theory of $\mathfrak{sl}_2[t]$.
\subsection{}We will determine the polynomial \eqref{poly} in two steps. In the first step we will reduce the computation of graded multiplicities in Demazure flags of fusion products to the computation of graded multiplicities in local Weyl modules only. To be more precise, we will show that there exists a polynomial in finitely many variables such that a suitable evaluation of that polynomial at elements of the form $\mathcal{V}^{1\rightarrow \xi_i}_{\bullet,\bullet}(q)$ for $1\leq i\leq \ell$ gives the desired polynomial $\mathcal{V}^{\xi\rightarrow m}_{n}(q)$. In the second step we give a combinatorial formula for graded multiplicities in Demazure flags of local Weyl modules in terms of Dyck paths. 
\begin{prop}\label{reductionpr} \mbox{}
\begin{enumerate}
\item Let $\xi_0\in \bn$ and $\xi=(\xi_1\geq \xi_2\geq \cdots \geq \xi_{\ell})$, $\ell\geq 2$, a partition such that $\xi_0\geq \xi_1$. Then  
$$\mathcal{V}_{|\xi|-2k}^{\xi\rightarrow \xi_0}(q)=\sum_{ 0=p_0\leq p_1\leq \cdots \leq p_{\ell-2}\leq p_{\ell-1}=k} q^{p_1+\cdots+p_{\ell-2}}\ \prod_{i=1}^{\ell-1} \mathcal{V}_{|\xi|_i-2p_{\ell-i-1},|\xi|_i-2p_{\ell-i}}^{\xi_{i}\rightarrow \xi_{i-1}}(q).$$
\item If $m,m'\in \bn$ and $s,n\in\bz_+$ are such that $m \geq m'$, then
$$\mathcal{V}^{m'\rightarrow m}_{s,n}(q)=\sum_{j=0}^{s}\left( \sum_{\ell=0}^j(-1)^{\ell+1}\hspace{-0,5cm} \sum_{0=p_0<p_1<\cdots<p_{\ell}<p_{\ell+1}=j}\ \prod_{i=0}^{\ell}\ \mathcal{V}^{1\rightarrow m'}_{s-2p_i,s-2p_{i+1}}(q)\right)\  \mathcal{V}^{1\rightarrow m}_{s-2j,n}(q).$$
\end{enumerate}
\proof
We prove the first part by an induction on $\ell$. If $\ell=2$, the discussion in Section \ref{section23} implies that $V(\xi)\cong D(\xi_1,|\xi|)$. Hence
$\mathcal{V}_{|\xi|-2k}^{\xi\rightarrow \xi_0}(q)=\mathcal{V}_{|\xi|,|\xi|-2k}^{\xi_{1}\rightarrow \xi_0}(q),$ and the induction begins. If $\ell>2$, we have
\begin{align*}
&& \mathcal{V}_{|\xi|-2k}^{\xi\rightarrow \xi_0}(q)&=\sum_{p=0}^k \mathcal{V}_{|\xi|-2p}^{\xi\rightarrow \xi_1}(q)\ \mathcal{V}_{|\xi|-2p,|\xi|-2k}^{\xi_1\rightarrow \xi_0}(q)&& \text{by \eqref{mat}}\\
&& &=\sum_{p=0}^k q^p\hspace{0.05cm}\mathcal{V}_{|\xi|_2-2p}^{\xi'\rightarrow \xi_1}(q)\ \mathcal{V}_{|\xi|-2p,|\xi|-2k}^{\xi_1\rightarrow \xi_0}(q)&& \text{by Lemma~\ref{ausc}}.
\end{align*}
The claim now follows by a trivial application of the induction hypothesis.

Now we prove the second part of the proposition. If $s-n\notin 2\bz_+$ both sides of the equation are zero and the statment is trivial. Otherwise assume that $s=n+2k$ for some $k\in\bz_+$. We can rewrite \eqref{mat} as a system of linear equations $Uv=b$, where $U=(\mathcal{V}_{i,j})_{0\leq i,j\leq k}$ is an upper triangular unipotent matrix of size $(k+1)\times (k+1)$ and
$$\mathcal{V}_{i,j}=\mathcal{V}^{1\rightarrow m'}_{s-2i,s-2j}(q),\ \ b_i=\mathcal{V}^{1\rightarrow m}_{s-2i,n}(q),\ \ v_i=\mathcal{V}^{m'\rightarrow m}_{s-2i,n}(q).$$
Now it is straightforward to check that the inverse of $U$ is again upper triangular unipotent whose entries are given by $f_{i,j}(\mathcal{V}),\  0\leq i,j\leq k,$
where $f_{i,j}(\mathcal{V})=\delta_{i,j}$ if $i\geq j$ and otherwise
$$f_{i,j}(\mathcal{V})=\sum_{r\geq 0}(-1)^{r+1} \sum_{i=p_0<p_1<\cdots<p_{r}<p_{r+1}=j}\mathcal{V}_{p_0,p_1}\mathcal{V}_{p_1,p_2}\cdots \mathcal{V}_{p_r,p_{r+1}}.$$
Now the first entry of the vector $v=U^{-1}b$ gives the desired result.
\endproof
\end{prop}
\begin{example} Let $V(\xi)$ be a truncated Wey module, i.e. $\xi=((d+1)^j,d^{N-j})$ where $N,d,j\in \bz_+$ are such that $0<j<N$ and $|\xi|=dN+j$. Applying Proposition~\ref{reductionpr} to this setting gives
$$\mathcal{V}_{|\xi|-2k}^{\xi\rightarrow \xi_0}(q)=\sum_{p=0}^k q^p\ \mathcal{V}_{(N-j)d,(N-j)d-2p}^{d \rightarrow d+1}(q)\ \mathcal{V}_{|\xi|-2p,|\xi|-2k}^{d+1 \rightarrow \xi_0}(q).
$$
\end{example}
\subsection{}After the reduction in Proposition \ref{reductionpr}, our second result focussed on a combinatorial formula for the graded multiplicities in a level $m$--Demazure flag of a local Weyl module. Recall the definition of admissible Dyck paths from Definition~\ref{maindef}.
\begin{thm}\label{mainthm1} Let $m\in\bn$, $n\in\bz_+$. We have, 
$$\mathcal{A}^{1\rightarrow m}_n(x,q)=\sum_{P\in \mathcal{D}_{m,n}}q^{\text{comaj}(P)}\ x^{d(P)},$$
where $d(P)$ denotes the number of down--steps of $P$. In particular, for any $s\in\bz_+$ we get
$$\mathcal{V}^{1\rightarrow m}_{s,n}(q)=\sum_{P\in\mathcal{D}_{m,n}(s)}q^{\text{comaj}(P)}.$$
\end{thm}
The proof of the above theorem will be postponed to Section~\ref{section5}, but to avoid technical difficulties in the rest of the paper we will handle the case $m=1$ seperately. If $m=1$ we have $A(1,n)=\{(0,1),(1,2),\dots,(n-1,n)\}$ and hence 
$$\mathcal{D}_{1,n}=\{\underbrace{11\cdots 1}_{n}\}=\mathcal{D}_{1,n}(n).$$ This yields,
$$\sum_{P\in \mathcal{D}_{1,n}(s)}q^{\text{comaj}(P)}=\delta_{s,n}$$
and the theorem follows immediately from \eqref{initial12}. \textit{So we will assume from now on that $m\geq 2.$} 
\begin{example}
If $n=3$, $m=2$ and $s=5$ we have $A(2,3)=\{(0,2),(2,3)\}$ and hence $\mathcal{D}_{2,3}(5)=\{10111, 11011\}$. We get $\mathcal{V}_{5,3}^{1\rightarrow 2}(q)=q^3+q^4$.
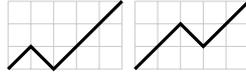
\begin{figure}[H]
\begin{tikzpicture}[scale=0.3]
\hspace{-1cm}
\draw [step=1,thin,gray!40] (0,0) grid (5,3);
\draw [very thick] (0,0) -- (1,1)--(2,0) --(5,3);
\end{tikzpicture}
\begin{tikzpicture}[scale=0.3]
\draw [step=1,thin,gray!40] (0,0) grid (5,3);
\draw [very thick] (0,0) -- (2,2) -- (3,1)--(5,3);
\end{tikzpicture}
\label{fig3}
\end{figure}
\end{example}
\begin{rem}We shall prove in Section~\ref{section5} a slightly more general version of Theorem~\ref{mainthm1}, namely we will give a combinatorial formula for the polynomials \eqref{poly} in terms of Dyck paths for every hook partition $\xi$ (see Theorem \ref{hook}).
\end{rem}
\subsection{} Before we proceed to the proof of Theorem~\ref{mainthm1} we will discuss several consequences. For simple finite-dimensional Lie algebras of non--simply laced type local Weyl modules are in general not isomorphic to level one Demazure modules. Nevertheless, Naoi proved in \cite[Sections 4 and 9]{Na11} that a local Weyl module (see for example \cite[Definition 3.1]{Na11} for a precise definition) admits a level one Demazure flag. If $\lie g$ is of type $B_n$ or $G_2$ one can use the $\mathfrak{sl_2}$--theory (in particular, the combinatorial formula in Theorem \ref{mainthm1}) to give combinatorial formulas for $\sum_p [W(\lambda): \tau_p^* D(1,\mu)]q^p$, where $\lambda,\mu$ are dominant integral $\lie g$--weights. This fact follows from an inspection of the proof given in \cite{Na11} and a precise statement can be found in \cite[Proposition 2.5]{CSSW14}.
\begin{cor}Let $\lie g$ of type $B_n$ and $m=2$ or $G_2$ and $m=3$. Further denote by $\alpha$ the unique simple short root with coroot $h_{\alpha}$. Let $\lambda,\mu$ two dominant integral weights of $\lie g$. We have
$$ [W(\lambda): \tau_p^* D(1,\mu)]= \begin{cases}|\{P\in \mathcal{D}_{m,\mu(h_{\alpha})}(\lambda(h_{\alpha})): \text{comaj}(P)=p\}|,& \lambda-\mu\in \bz_+\alpha\\
0,& \text{otherwise.}
\end{cases}$$
\hfill\qed
\end{cor}
\subsection{} In the special case when $m=2$, we can derive another combinatorial description of the graded multiplicities in terms of bounded partitions which we record in the next lemma. This result can also be obtained by a straightforward calculation combining \cite[Theorem 3.3]{CSSW14} and \cite[Chapter 3]{And98} as pointed out in \cite[Section 4.3]{MJ16}. Our proof is different and uses only the combinatorial formula stated in Theorem \ref{mainthm1}. For integers $a,b,c\in\bz$, denote by $\rho^{a}_{b}(c)$ the set of partitions of $c$ with at most $a$ parts such that each part is bounded by $b$. Recall the integers $n_0,n_1$ from Section~\ref{section23}
\begin{lem}\label{lpart} Let $k,n\in \bz_+$. Then we have
$$\mathcal{V}^{1\rightarrow 2}_{n+2k,n}(q)=\sum^{k(k+n)}_{\ell=0}   |\rho^{k}_{n_1}(k(k+n)-\ell)|\hspace{0,04cm} q^{\ell}.$$
\end{lem}
\proof If $0\leq n\leq 1$ we have $n_1=0$ and hence the right hand side is equal to $q^{k^2}$. Since each path in $\mathcal{D}_{n}^{1}(n+2k)$ has height $\leq 1$ the same is true for the left hand side and the lemma is immediate; so assume that $n>1$. Let $P\in \mathcal{D}_{2,n}(n+2k)$ such that $\text{comaj}(P)=\ell$. Since 
$$A(2,n)=\{(i,i+2),(p,p+1): 0\leq i\leq n_1-1,\ n_1<p\leq n-1\}$$
we see that any admissible path has no two adjacent down--steps and the last $n_1$ steps are up--steps. Hence $P$ is uniquely determined by the $x$--coordinates of all its peaks and the number of down--steps is equal to the number of peaks (= $k$). Assume that $i_1,\dots,i_{k}$ are these $x$--coordinates, then we must have  $i_{p+1}\geq i_p+2$ for all $1\leq p<k$ and
$$1\leq i_1<i_2<\cdots<i_{k}<2k+n_1,\ \ \sum_{p=1}^{k} i_p=(n+2k)k-\ell.$$
Now substituting $i_p'=i_p-2p+1$ we get a partition 
$$0\leq i'_1\leq i'_2\leq \cdots\leq i'_{k}\leq n_1, \ \ \sum_{p=1}^{k} i'_p=k(k+n)-\ell,$$
which obviously gives a bijective correspondence $$\{P\in \mathcal{D}_{2,n}(n+2k) : \text{comaj}(P)=\ell\}\cong \rho^{k}_{n_1}(k(k+n)-\ell).$$
\endproof

\subsection{}Our next corollary answers the following question. What is the set of all powers $q$ which appear with nonzero coefficients in the polynomial $\mathcal{V}^{1\rightarrow m}_{s,n}(q)$? We first determine the degree of $\mathcal{V}^{1\rightarrow m}_{s,n}(q)$ as we did in Lemma \ref{lpart} for $m=2$.  By \eqref{initial12} we can assume that $s=n+2k$ for some $k\in \bz_+$. The path 
\begin{equation}\label{4}P=\underbrace{1010\cdots 1010}_{2k}\underbrace{111\cdots 11}_{n}\end{equation}
with exactly $k$ up--down steps at the beginning and $n$ up--steps at the end is contained in $\mathcal{D}_n(s)$. Note that all peaks of $P$ are at height $1$ and $(0,1)\notin A(m,n)$, which implies that $P$ is admissible. Moreover,
\begin{equation}\label{2}\text{comaj}(P)=\sum_{i=1}^k(s-2i+1)=k(n+k).\end{equation}
If $P'\in \mathcal{D}_n(s)$ is another Dyck path with peaks at length $s_1,\dots,s_p$, we have $s_i\geq 2i-1,\ i=1,\dots,p$ and $p\leq k$. Thus
\begin{equation}\label{3}\text{comaj}(P)=\sum_{i=1}^k(s-2i+1)\geq \sum_{i=1}^p(s-s_i)+\sum_{i=p+1}^k(s-2i+1)\geq\text{comaj}(P')\end{equation}
with equality if and only if $P=P'$.
\begin{cor}\label{conint}We have that $\mathcal{V}^{1\rightarrow m}_{s,n}(q)=0$ if $s-n\notin 2\bz_+$ and otherwise it is a monic polynomial of degree $(s-n)(s+n)/4$. Moreover, the powers of $q$ which appear in the polynomial $\mathcal{V}^{1\rightarrow m}_{s,n}(q)$ form an interval of consecutive integers iff one of the following conditions hold
(i) $n\neq 0$, \ \  (ii) $m=2$, \ \ (iii) $s=n$.
\end{cor}
\proof
The first part of the corollary follows from the discussion preceding the corollary, see \eqref{2} and \eqref{3}. We start by proving the backward direction of the second part. Let $n\neq 0$ and $\tilde{P}\neq P$ an arbitrary path in $\mathcal{D}_{m,n}(s)$, where $P$ is as in \eqref{4}. We shall show that $q^{\text{comaj}(\tilde{P})+1}$ appears with a nonzero coefficient in the polynomial $\mathcal{V}^{1\rightarrow m}_{s,n}(q)$. 
Since $\tilde{P}\neq P$, we must have two consecutive up-steps in $\tilde{P}$ followed by at least one down-step. Thus $\tilde{P}$ is of the form 
$$\tilde{P}=\cdots \underbrace{11\cdots 11}_{\ell}\underbrace{00\cdots 00}_{r}1\cdots,\ \ell\geq 2,\ \ r\geq 1,$$
and we can assume without loss of generality that the path preceding the $\ell$--up steps is contained in $\mathcal{D}_0$, i.e. it ends at the $x$--axis.
\textit{Case 1:} Suppose that $r<\ell$. We set $$\bar{P}=\cdots \underbrace{11\cdots 11}_{\ell-1}\underbrace{00\cdots 00}_{r}11\cdots$$
and obtain $\bar{P}\in \mathcal{D}_n^N(s)$ and $\text{comaj}(\bar{P})=\text{comaj}(\tilde{P})+1.$
\begin{figure}[H]
\begin{tikzpicture}[scale=0.3]
\draw [step=1,thin,gray!40] (0,0) grid (19,7);
\draw [very thick] (0,0) -- (3,3) -- (6,0) -- (10,4)--(13,1)--(17,5)--(19,7);
\draw [dotted] (9,3) -- (12,0) -- (13,1);
\node [left=1pt] at (8,2) {\scriptsize $\ell$};
\node [left=1pt] at (12,2) {\scriptsize $r$};
\end{tikzpicture}
\caption{Solid path $\tilde{P}$ and dotted path $\bar{P}$}
\end{figure}

To complete the proof in this case, we have to show that $\bar{P}$ is admissible. Again, there is nothing to show if $n<m$; so suppose that $n\geq m$. We denote by $z$ the peak of $\tilde{P}$ after the $\ell$ up--steps and set for simplicity $b:=\text{ht}_{z}(\tilde{P})$. It is clear that $\bar{P}$ has a peak at height $b-1$. If $b\leq m$, there is nothing to show. Otherwise we must show that $\text{ht}_{\tilde z}(\bar{P})>\tilde a$ for all points $\tilde z\in \text{supp}(\bar{P})$ with $\ell_{\tilde z}(\bar{P})\geq \ell_{z}(\tilde{P})-1$,
where $\tilde a$ is the unique non--negative integer such that $(\tilde a,b-1)\in A(m,n)$. For a point $\tilde z\in \text{supp}(\bar{P})$ we have that $\tilde z\in \text{supp}(\tilde{P})$ or $(\ell_{\tilde z}(\bar{P})+1,\text{ht}_{\tilde z}(\bar{P})+1)\in \text{supp}(\tilde{P})$. In either case,
$$\text{ht}_{\tilde z}(\bar{P})> a-1\geq \tilde a,$$
where first inequality follows by our assumption $\tilde{P}\in \mathcal{D}_{m,n}(s)$ and the second inequality follows from \eqref{1}.

\textit{Case 2:} In this case we can suppose that each peak in $\tilde{P}$ returns to the $x$--axis. Since $\tilde{P}\neq P$, we have that $\tilde{P}$ contains at least two down--steps. In particular $m\geq 3$, since $(0,m)\in A(m,n)$ or when $n<m$ the height of $\tilde{P}$ is at most $m-1$. Hence 
$$\tilde{P}=\cdots 1\underbrace{00\cdots 00}_{r}\underbrace{1010\cdots 1010}_{t}\underbrace{11\cdots 11}_{n},\ \ r\geq 2,\ \ t\in 2\bz_+.$$
If $t\neq 0$, we can set 
$$\tilde{P}=\cdots 1\underbrace{00\cdots 00}_{r-1}100\underbrace{10\cdots 1010}_{t-2}\underbrace{11\cdots 11}_{n},$$
which is admissible and has the desired property. If $t=0$, we set
$$\bar{P}=\cdots 1\underbrace{11\cdots 11}_{N_{0,1}}\underbrace{00\cdots 00}_{r-1}\underbrace{11\cdots 11}_{n-2N_{0,1}}0\underbrace{11\cdots 11}_{N_{0,1}},$$
where $N_{0,1}=n_1+\delta_{n_0,m-1}$. We get
$$\text{comaj}(\bar{P})=\text{comaj}(\tilde{P})-(n+r)+(n+r-N_{0,1})+(N_{0,1}+1)=\text{comaj}(\tilde{P})+1.$$
\begin{center}
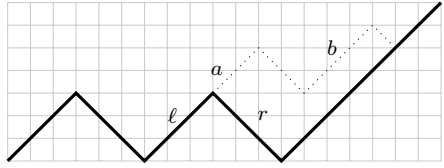
\begin{figure}[H]
\begin{tikzpicture}[scale=0.3]
\draw [step=1,thin,gray!40] (0,0) grid (19,7);
\draw [very thick] (0,0) -- (3,3) -- (6,0) -- (9,3)--(12,0)--(19,7);
\draw [dotted] (9,3) -- (11,5) -- (13,3)--(16,6)--(17,5);
\node [left=1pt] at (8,2) {\scriptsize $\ell$};
\node [left=1pt] at (12,2) {\scriptsize $r$};
\node [left=1pt] at (10,4) {\scriptsize $a$};
\node [left=1pt] at (15,5) {\scriptsize $b$};
\end{tikzpicture}
\caption{Solid path $\tilde{P}$ and dotted path $\bar{P}$ with $a=N_{0,1}$, $b=n-2N_{0,1}$}
\end{figure}
\end{center}
Again, we must prove that $\bar{P}$ is admissible. Recall that $r<m$, since $N=m-1$ or $(0,m)\in A(n,m)$. There are two new peaks $z_1$ and $z_2$ in the path $\bar{P}$, where 
$$z_1=(s-n-r+N_{0,1},r+N_{0,1}),\ \ z_2=(n-N_{0,1}-1,n-N_{0,1}+1).$$ We know that 
$$(n-N_{0,1}-1,n-N_{0,1}+1)\in A(n,m)$$ and hence the peak $z_2$ does not violate property \eqref{admissible}. If $r+N_{0,1}<m$, there is nothing to check for $z_1$ and otherwise we have $(r+N_{0,1}+\delta_{n_0,m-1}-m,r+N_{0,1})\in A(m,n)$ and the claim follows with
$$\text{ht}_{z'}(\bar{P})\geq N_{0,1}+1>r+N_{0,1}+\delta_{n_0,m-1}-m,\ \ \text{ $\forall z'$ with }\ell_{z'}(\bar{P})\geq s-n-r+N_{0,1}.$$
Hence in the case when $n\neq 0$ we have that $\mathcal{V}^{1\rightarrow m}_{s,n}(q)$ forms an interval of consecutive integers.  If $n=0$ and $m=2$, the statement is clear since $\mathcal{D}_{2,0}(s)=\{P\}$.  Similarly there is nothing to prove if $n=0=s$ and hence the backward direction is established.

Now we show that $\mathcal{V}^{1\rightarrow m}_{s,n}(q)$ does not form an interval of consecutive integers provided that $n=0$, $m\geq 3$ and $s\neq n$. For this we consider a path of minimal degree in $\mathcal{D}_{m,n}(s)$ given by 
$$Q=\underbrace{11\cdots 11}_{m-1}\underbrace{00\cdots 00}_{m-1}\cdots\underbrace{11\cdots 11}_{m-1}\underbrace{00\cdots 00}_{m-1}\underbrace{11\cdots 11}_{p}\underbrace{00\cdots 00}_{p},$$
where $(s-n)/2=\ell (m-1)+p$,\  $0\leq p<m-1$. Then it is easy to show that there is no path in $\mathcal{D}_{m,n}(s)$ whose comajor is given by $\text{comaj}(Q)+1$.
\endproof
The representation theoretic meaning of the above corollary is quite surprising.
\begin{lem} The maximal degree in which a level $m$--Demazure module appears in a level $m$--Demazure flag of a local Weyl module is independent of $m$ and depends only on its highest weight. Moreover, given non--negative integers $s,p_1,p_2\in \bz_+$ with $p_1<p_2$ and $m,n\in \bn$, we have an equivalence
$$[W_{\text{loc}}(s):\tau_p^{*}D(m,n)]\neq 0\ \ \forall p\in [p_1,p_2]\cap \bz\iff [W_{\text{loc}}(s):\tau_{p}^{*}D(m,n)]\neq 0\ \ \forall p\in \{p_1,p_2\}.$$
\hfill\qed
\end{lem}
\subsection{}\label{section33}Under the restriction $m\geq s$ the module $D(m,s)$ is irreducible as an $\mathfrak{sl}_2$--representation (see for example \cite{CSSW14}). Hence the level $m$--Demazure flag becomes the usual Jordan--H\"older series. In this situation graded multiplicities in fusion products were studied in 
\cite{Ke04} and are described by co--charge Kostka--Foulkes polynomials.
Let $\langle\cdot,\cdot\rangle$ be the Hall inner product on $\Lambda$,
the ring of symmetric functions. The modified Hall--Littlewood polynomials 
$Q'_{\mu}=Q'_{\mu}(x;q)$ are defined as the basis of $\Lambda[q]$ 
dual to the ordinary Hall--Littlewood polynomials 
$P_{\lambda}=P_{\lambda}(x;q)$ with respect to
$\langle\cdot,\cdot\rangle$:
\[
\langle P_{\lambda},Q'_{\mu}\rangle=\delta_{\lambda,\mu}.
\]
The Kostka--Foulkes polynomials $K_{\lambda,\mu}(q)$ are then defined as 
\[
Q'_{\mu}(x;q)=\sum_{\lambda} K_{\lambda,\mu}(q)s_{\lambda}(x),
\]
where the $s_{\lambda}(x)$ are the Schur functions.
The closely related co--charge Kostka polynomials $\tilde{K}_{\lambda,\mu}(q)$
are given by
\[
\tilde{K}_{\lambda,\mu}(q)=q^{n(\mu)} K_{\lambda,\mu}(1/q),
\]
where $n(\mu):=\sum_{i\geq 1} (i-1)\mu_i$. The following result goes back to \cite{Ke04} and an alternative proof can be given by using $q$--Bernstein operators (see \cite{Jin91}) to show that the co--charge Kostka polynomial satisfies a recursion similar to \cite[Theorem 5]{CV13}. We thank Ole Warnaar for helping us with this observation.
\begin{lem}\label{ked}Let $\xi$ a partition and $m\in \bn$ sucht that $m\geq |\xi|$. Set $\mu=(\frac{|\xi|+n}{2}, \frac{|\xi|-n}{2} )$. We have 
 $$\mathcal{V}_{n}^{\xi\rightarrow m}(q)=\tilde{K}_{\mu,\xi}(q).$$
 \hfill\qed
\end{lem}
In the rest of this section we discuss an interesting corollary. Let $s=n+2k$ for some $k\in\bz_+$. We claim that 
\begin{equation}\label{6}\mathcal{D}_{m,n}(s)=\mathcal{D}_n(s),\ \text{ if $m\geq s$}.\end{equation}
If $k=0$, the claim is immediate. Otherwise, note that any path $P\in \mathcal{D}_n(s)$ has height at most $n+k$ and $n+k\leq m-1$. Hence $\mathcal{D}_n(s)=\mathcal{D}^{m-1}_n(s)= \mathcal{D}_{m,n}(s)$, since $N=m-1$ and $A(m,n)=\emptyset$. This shows \eqref{6}. There is a bijection between $\mathcal{D}_n(s)$ and the set of standard Young tableaux $SYT(\lambda)$ of shape $\lambda=(n+k,k)$ as follows:
\begin{equation}\label{bijtoSYT}P=a_1\cdots a_{s}\mapsto T(P),\end{equation}
where we put $i$ into the first row if $a_i=1$ and otherwise into the second row. See Figure~\ref{fig2} for an example.
\begin{figure}[H]
\begin{tikzpicture}[scale=0.3]
\draw [step=1,thin,gray!40] (0,0) grid (7,2);
\draw [very thick] (0,0) -- (2,2) -- (3,1) -- (4,2)--(6,0)--(7,1);
\end{tikzpicture}
$\mapsto \young(1247,356)$
\caption{}
\label{fig2} 
\end{figure}
As mentioned in Remark~\ref{rem1} there is also a comajor statistics on standard Young tableaux which is often useful. Given a tableaux $T$ of shape $\lambda$, a descent of $T$ is a value $i$, $1\leq i\leq s$, for which $i+1$ occurs in one of the rows below $i$. Define $\text{comaj}(T)=\sum_i s-i$, where the sum runs over the descents of $T$.  Stanley proved \cite[pg. 363]{S99}, 
\begin{equation}\label{stan}s_{\lambda}(1,q,q^2,\dots)=\frac{1}{(q;q)_{s}}\sum_{T\in SYT(\lambda)}q^{\text{comaj}(T)}.\end{equation}
\begin{cor}\label{princp}Let $\lambda=(n+k,k)$, $\xi=(1^{n+2k})$ and $\mu=(\frac{s+n}{2}, \frac{s-n}{2} )$. Then
$$\tilde{K}_{\mu,\xi}(q)=(q;q)_{n+2k} \ s_{\lambda}(1,q,q^2,\dots).$$
\proof Let $V_{\mathfrak{sl}_2}(n)$ the $(n+1)$--dimensional irreducible $\mathfrak{sl}_2$--representation. Theorem~\ref{mainthm1} together with \eqref{6} imply
\begin{equation}\label{7}\sum_{p\geq 0}[W_{\text{loc}}(n+2k): \tau_p^* V_{\mathfrak{sl}_2}(n)]\hspace{0,03cm}  q^p=\sum_{P\in \mathcal{D}_n(s)}q^{\text{comaj}(P)}.\end{equation}
Since \eqref{bijtoSYT} preserves the comajor statistics we get that \eqref{7} is equal to $\sum_{T\in SYT(\lambda)}q^{\text{comaj}(T)}$ and the claim follows with \eqref{stan} and Lemma~\ref{ked}.
\endproof
\end{cor}

\subsection{}In his last letter to G.H. Hardy, S. Ramanujan listed 17 functions which he called mock theta functions \cite{Ram88}. There are some number theoretic interpretations of some of these functions in the literature. For example, one of the third order mock theta functions has been interpreted as the generating function for partitions into odd parts without gaps \cite{Fi88}.  Agarwal gave an interpretation of some fifth order mock theta functions in terms of $n$--color partitions \cite{Ag04} and later in terms of lattice paths \cite{Ag05}. Thanks to \cite[Theorem 1.6]{BCSV15} we have a relationship between certain specializations of the series $\mathcal{A}_n^{1\rightarrow 3}(x,q)$, $n\in\{0,1,2\}$ and the fifth order mock theta functions $\phi_0(q),\phi_1(q),\psi_0(q),\psi_0(q)$, where
\begin{align*}
\phi_0(q) = \sum_{n=0}^\infty q^{n^2} (-q;q^2)_n,\ \
\phi_1(q) = \sum_{n=0}^\infty q^{(n+1)^2} (-q;q^2)_n,\\
\psi_0(q) = \sum_{n=0}^\infty q^{\frac{(n+1)(n+2)}{2}} (-q;q)_n,\ \
\psi_1(q) = \sum_{n=0}^\infty q^{\frac{n(n+1)}{2}}(-q;q)_n.
\end{align*}
As an immediate consequence (combining \cite[Theorem 1.6]{BCSV15} and Theorem~ \ref{mainthm1}) we get a new interpretation of their coefficients. We emphasize that our interpretation is different from \cite{Ag05}.
\begin{cor}\label{mockcoef}For $n\in\bz_+$ let $res_2(n)$ the remainder of $n$ modulo 2. We have
\begin{enumerate}
\item $\psi_1(q)=\sum_{n=0}^{\infty}|\{P\in \mathcal{D}^2_{1}: \text{comaj}(P)=n\}|q^n$,
\vspace{0,2cm}
\item $\psi_0(q)=\sum_{n=0}^{\infty}|\{P\in \mathcal{D}^2_{1}: \text{comaj}(P)+\lceil\frac{\ell(P)}{2}\rceil=n\}|q^n,$
\vspace{0,2cm}

\item $\phi_0(q)=\sum_{n=0}^{\infty}|\{P\in \mathcal{D}^2_{2\text{res}_2(n)}: \text{comaj}(P)=\lfloor n \rfloor\}|q^n$
 \vspace{0,2cm}

\item $\phi_1(q)=\sum_{n=0}^{\infty}|\{P\in \mathcal{D}^2_{2(1-\text{res}_2(n))}: \text{comaj}(P)+\lceil\frac{\ell(P)}{2}\rceil+(1-\text{res}_2(n))=\lfloor n \rfloor\}|q^n.$
\end{enumerate}
\hfill\qed
\end{cor}
\subsection{} The connection to mock theta functions \cite[Theorem 1.6]{BCSV15} after specializing the generating series $\mathcal{A}_n^{1\rightarrow 3}(x,q)$, $n\in\{0,1,2\}$ is quite surprising. We emphasize that the relationship could only be made for one reason: the series $\mathcal{A}_n^{1\rightarrow 3}(x,q)$, $n\in\{0,1,2\}$ can be expressed in terms of $q$--binomial coefficients \cite[Sections 1.4 and 1.6]{BCSV15}. In this section we want to generalize the situation and give a formulae for $\mathcal{A}_n^{1\rightarrow m}(x,q)$, $n\in\{0,\dots,m-1\}$ in terms of $q$--binomials, but the connection to number theory will be discussed elsewhere. The key ingredient in the proof of the next corollary is Theorem ~\ref{mainthm1} and a result of Krattenthaler and Mohanty \cite[Theorem 1]{KM93}.
\begin{cor}\label{KrM}Let $m\in \bn$ and $n\in \bz_+$ such that $n<m$. We have that $\mathcal{V}_{n+2k,n}^{1\rightarrow m}(q)$ equals
\begin{equation*}\label{b}\sum_{\ell=0}^k\sum_{s=-\ell}^{\ell} q^{\ell^2+s(sm+n+1)} \text{det} \begin{pmatrix}
\qbinom{k-s(m-1)}{\ell+s}_q & 
\qbinom{k+s(m-1)-1}{\ell-s}_q  \\
q^{-2s(n+1)}\qbinom{n+1+k-s(m-1)}{\ell+s}_q & \qbinom{n+k+s(m-1)}{\ell-s}_q \\
\end{pmatrix}.\end{equation*}
\end{cor}
\proof Since $n<m$ we have $\mathcal{D}_{m,n}(n+2k)=\mathcal{D}^{m-1}_{n}(n+2k)$ and hence Theorem ~\ref{mainthm1} implies
$$\mathcal{V}_{n+2k,n}^{1\rightarrow m}(q)=\sum_{P\in \mathcal{D}^{m-1}_{n}(n+2k)} q^{comaj(P)}.$$
It means that we have to count the set of Dyck paths from $(0,0)$ to $(n+2k,n)$ of height $< m$ by its comajor index. A special case of Theorem 1 of \cite{KM93} counts the set of Dyck paths from $(0,0)$ to $(n+2k,n)$ of height $<m$ with respect to a different statistics, namely by its major and descent index. The descent index $\text{des}(P)$ counts the number of peaks of $P$ and hence
\begin{equation}\label{desmaj}\text{maj}(P)=\text{des}(P)(n+2k)-\text{comaj}(P).\end{equation}
In view of \eqref{desmaj} we only have to substitute the formula in Theorem ~\ref{mainthm1} by $q\mapsto 1/q$, $x\mapsto x^{n+2k}$ and finally $x=q$.
\endproof
\begin{rem}
The advantage of Corollary ~\ref{KrM} is that it gives some hints to compute graded multiplicities for higher rank Lie algebras. As remarked in the introduction level $m'$--Demazure modules for $\mathfrak{sl}_n$ also admit a level $m$--Demazure flag if $m\geq m'$. So the polynomials \eqref{poly} can be defined in a similar manner; also the corresponding generating series. We do not know yet whether the size of the matrix in Corollary ~\ref{KrM} has something to do with the rank of the underlying Lie algebra or the fact that our model is based on certain two dimensional lattice paths. One possibility for higher rank Lie algebras would be to study multidimensional lattice paths (depending on the rank of the Lie algebra), see for example \cite{HM79}. We will address this in a future publication.  
\end{rem}


\section{Proof of Theorem~\ref{mainthm1}}\label{section5}
In this section we prove our main theorem. 
\subsection{} We first show that our combinatorial model gives the correct numerical multiplicity and a proof of that can be given purely combinatorial, which we will demonstrate in this subsection. The generating function $\mathcal{G}_{m,n}^N(x):=\sum_{s\geq 0}|\mathcal{D}_{m,n}(s)| x^{(s-n)/2}$ is again a rational function which can be expressed in terms of Chebyshev polynomials. The next lemma can be viewed as a generalization of Lemma \ref{Ges}.
\begin{prop}\label{numerical} We have
$$\mathcal{G}_{m,n}^N(x)=\frac{p_{m-n_0-1}(x)}{p_m(x)^{n_1+1}}=\mathcal{A}_n^m(x,1).$$
In particular,
$$\mathcal{V}_{s,n}^{1\rightarrow m}(1)=|\mathcal{D}_{m,n}(s)|.$$
\end{prop}
\proof
If $n<m$, the statement follows from Lemma~\ref{Ges}; so assume that $n\geq m$. Let $(i_1,\dots,i_{n-m+1})$ be the integers from Definition~\ref{maindef} and set $i_0=0$, $i_{n-m+2}=n$ and $\tilde{\mathcal{G}}_{m,n}^N(x)=x^n\mathcal{G}_{m,n}^N(x)$. Let $P$ be a path contained in $\mathcal{D}_{m,n}$. We want to factor $P$ as 
\begin{equation}\label{factorpath}P=P_11P_21\cdots P_{n-m+1}1 P_{n-m+2},\end{equation} where $P_{\ell}$ is a path of height at most $m-i_{\ell-1}+(\ell-2)+(\delta_{\ell,1}-1)$ that ends at height $i_{\ell}-i_{\ell-1}+(\delta_{\ell,1}-1)$. This would give, together with Lemma~\ref{Ges},
$$\tilde{\mathcal{G}}_{m,n}^N(x)=\prod^{n-m+2}_{\ell=1}x^{i_{\ell}-i_{\ell-1}}\frac{p_{m-i_{\ell}+(\ell-2)}(x)}{p_{m-i_{\ell-1}+(\ell-2)+\delta_{\ell,1}}(x)}=\frac{p_{m-n_0-1}(x)}{p_m(x)^{n_1+1}}.$$
The rest of the claim follows from \cite[Corollary 1.3]{BCSV15}. 
Now we make precise how we factor the path as in \eqref{factorpath}. The path $P$ starts at the origin and we follow the path until a point $z\in \text{supp}(P)$ with $z=(\ell_z(P),i_1)$. Since $P$ ends at height $n$ and $i_1<n$, the existence of $z$ is clear. Among all points we choose $z$ such that $\ell_z(P)$ is maximal, i.e. we take the right most point of height $i_1$.  We define this part of the path by $P_1$ and note that the height of $P_1$ is at most $m-1$, since $(i_1,m)\in A(m,n)$. By the maximality of $z$, we must have that $P_1$ is immediately followed by an up--step. Hence $P=P_11Q$, where $Q$ is the remaining part of $P$. We continue the process with $Q$, but note that this path is shifted up by $(i_1+1)$--units. Again we follow the path $Q$ until a point $z$ on the path with $z=(\ell_z(P),i_2)$ and $\ell_z(P)$ is maximal with this property (the right most point with height $i_2$.) This part of the path we define as $P_2$ and again the height of $P_2$ is at most $m-(i_1+1)$, since $(m+1,i_2)\in A(m,n)$ (remember the shift by $(i_1+1)$). Continuing this way gives the desired factorization.
\begin{figure}[H]
\begin{tikzpicture}[scale=0.3]
\draw [step=1,thin,gray!40] (0,0) grid (21,5);
\draw [very thick,blue] (0,0) -- (3,3) -- (5,1) -- (7,3)--(10,0);
\draw [very thick] (10,0) -- (11,1);
\draw [very thick,green](11,1)-- (14,4)--(17,1)--(18,2);
\draw [very thick] (18,2) -- (19,3);
\draw [very thick,red] (19,3) -- (21,5);
\end{tikzpicture}
\caption{A factorization $P={\color{blue}P_1} 1 {\color{green}P_1} 1 {\color{red}P_3}$ where $n_1=1, n_0=1, m=4$.}
\end{figure}
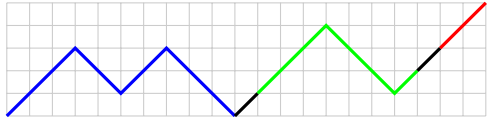
\endproof
\subsection{}
It is clear that each Dyck path in $\mathcal{D}_n(s)$ is uniquely determined by the coordinates of all its peaks. Let $P$ a path with $d$ peaks whose coordinates are
$$B_1=(x_1,y_1),\dots,B_d=(x_{d},y_{d}),\ \  x_1<x_2<\cdots< x_d.$$ 
We define inductively a subset $\eta(P)\subseteq\{1,\dots,d\}$ as follows. 
Set $i_1=1$ and define 
\begin{equation*}i_{t+1}:=\min\{i_t<u\leq d : y_u\geq y_{i_t}\},\ \ 1\leq t\leq \ell-1,\end{equation*}
where $\ell$ denotes the maximal non--negative integer such that $i_{\ell}$ exists. Define 
\begin{equation}\label{incsubsetofpeaks}
\eta(P):=\{i_1,\dots,i_{\ell}\}.\end{equation}
The meaning of the sequence $\eta(P)$ in the language of Dyck path is as follows: given a peak $B_{i_{t}}$ let $B_{i_{t+1}}$ the unique peak determined by the property that $B_{i_{t+1}}$ is the left most peak to the right of $B_{i_{t}}$ which is weakly above $B_{i_{t}}$. This sequence will be needed in the proof of Proposition \ref{eigenschaft1}.
\begin{example}
Let $P\in \mathcal{D}_7(20)$ the unique path with peaks
$$B_1=(3,3), B_2=(6,2), B_3=(10,4), B_4=(15,5), B_5=(17,5).$$
By definition we have $i_1=1$. Since $B_{2}$ is not weakly above $B_{i_1}$ but $B_{3}$ is weakly above $B_{i_1}$, we get $i_2=3$. Continuing in this way we get $\eta(P)=\{1,3,4,5\}$.
\begin{figure}[H]
\begin{tikzpicture}[scale=0.3]
\draw [step=1,thin,gray!40] (0,0) grid (20,7);
\draw [very thick] (0,0) -- (3,3) -- (5,1) -- (6,2)--(7,1)--(10,4)--(12,2)-- (15,5) -- (16,4)-- (17,5)--(18,4)--(19,5)--(20,6);
\end{tikzpicture}
\caption{$i_1=1$, $i_2=3$, $i_3=4$ and $i_4=5$}
\end{figure}
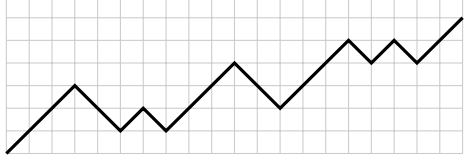
\end{example}
Similarly, a Dyck path $Q\in \mathcal{D}_n(s)$ is also uniquely determined by the coordinates of all its valleys. So assume that $Q$ has $d$ valleys whose coordinates are 
$$C_1=(x_1,y_1),\dots, C_d=(x_{d},y_{d}),\ \ x_1<x_2<\cdots< x_d.$$  
This time we set $j_1=d$ and define inductively a subset $\zeta(Q)\subseteq\{1,\dots, d\}$ as follows.  Set
\begin{equation*}j_{t+1}:=\max\{1\leq u< j_t : y_u\leq y_{j_t}\},\ \ v-1\leq t\leq d,\end{equation*}
where $v$ denotes the maximal non--negative integer such that $j_{v}$ exists.  Define 
\begin{equation}\label{incsubsetofvalleys}
\zeta(P):=\{j_1,\dots,j_{v}\}.\end{equation}
The meaning of the sequence $\zeta(P)$ in the language of Dyck paths is as follows: given a valley $C_{i_{t}}$ let $C_{i_{t+1}}$ the unique valley determined by the property that $C_{i_{t+1}}$ is the right most valley to the left of $C_{i_{t}}$ which is weakly below $C_{i_{t}}$. 
\subsection{}Let $_{j}\mathcal{D}_{m,n}(s)\subset \mathcal{D}_{m,n}(s)$ be the subset of Dyck paths which start with $j\in \bn$ number of up--steps, i.e.
$$_{j}\mathcal{D}_{m,n}(s)=\{a_1\cdots a_{s}\in\mathcal{D}_{m,n}(s): a_1=a_2=\cdots=a_{j}=1\}.$$ 
Furthermore, set 
$$_{j}\mathcal{E}_{s,n}^{1\rightarrow m}(q)=\sum_{P\in \hspace{0,02cm}  _{j}\mathcal{D}_{m,n}(s)} q^{\text{comaj(P)}}.$$
\begin{prop}\label{eigenschaft1}Let $n,s\geq m$ and $s-n\in 2\bz_+$. Then, there exists a bijection 
$$\Psi  : \mbox{$_m\mathcal{D}_{m,n}(s)$} \rightarrow \mathcal{D}_{m,n-m}(s-m)$$
 such that $\text{comaj}(P)=\text{comaj}(\Psi(P))+(s-n)/2$.
\proof
Set $s=n+2k$ for some $k\in \bz_+$. Let $P\in \mbox{$_m\mathcal{D}_{m,n}(s)$}$ and assume that $P$ has $d$ peaks whose coordinates are 
$$B_1=(x_1,y_1),\dots,B_d=(x_{d},y_{d}),\ \ x_1<x_2<\cdots< x_d.$$ 
Furthermore, let $\eta(P)=\{i_1,i_2,\dots,i_{\ell}\}$ as in \eqref{incsubsetofpeaks}.
Let $\Psi(P)$ the unique path in $\mathcal{D}_{m,n-m}(s-m)$ with $d$ peaks whose coordinates are given by
$$B^{\Psi}_i=(x_i+\ell_i-m,2y_{i_t}-y_i+\ell_i-m),\ \ 1\leq t\leq \ell,\ \ i_t\leq i < i_{t+1},$$
where $\ell_i$ denotes the number of down--steps right after peak $B_i$ and $i_{\ell+1}:=d+1$.

This map has a pictorial interpretation, namely we reflect the path $P$ piecewise. To be more precise, we draw $\ell$ horizontal lines connecting the points $(x_{i_t},y_{i_t})$ with $(x_{i_t}+s_{i_t},y_{i_t})$, $1\leq t\leq \ell$, where $s_{i_t}=\min\{s\in \bn : (x_{i_t}+s,y_{i_t})\in \text{supp}(P)\}$. Then we reflect the path locally at these lines and ignore the first $m$ up-steps (see Figure~\ref{fig6} for an example).
\begin{figure}[H]
\begin{tikzpicture}[scale=0.3]
\draw [step=1,thin,gray!40] (0,0) grid (12,9);
\draw [very thick] (0,0) -- (5,5) -- (6,4) -- (8,6)--(10,4)--(12,6);
\draw [green] (5,5) -- (7,5);
\draw [green] (8,6) -- (12,6);
\draw [dotted] (5,5) -- (6,6) -- (7,5)--(10,8)--(12,6);
\node [left=1pt] at (7,5) {\scriptsize $\ell_1$};
\node [left=1pt] at (10.5,5.5) {\scriptsize $\ell_2$};
\node [left=1pt] at (1.5,7.5) {\scriptsize $P$};
\end{tikzpicture}
\hspace{2cm}
\begin{tikzpicture}[scale=0.3]
\draw [step=1,thin,gray!40] (0,0) grid (7,5);
\draw [very thick] (0,0) -- (1,1) -- (2,0) -- (5,3)--(7,1);
\node [left=1pt] at (3,3.5) {\scriptsize $\bold{\Psi}(P)$};
\end{tikzpicture}
\caption{The case $m=5$, $k=3$, $n=6$. The path is reflected at the green lines.}
\label{fig6}
\end{figure}
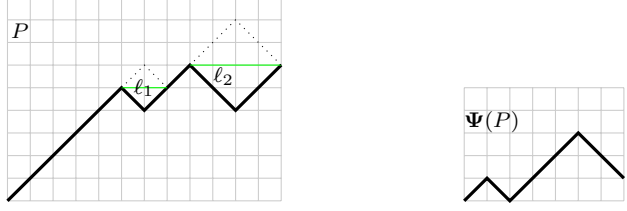
Since $\ell_1+\cdots+\ell_d=k$, we immediately get $\text{comaj}(\Psi(P))=\text{comaj}(P)-k$. In order to show that $\Psi$ is well--defined we note that $(x_i+\ell_i,y_i-\ell_i)\in \text{supp}(P).$ Fix $t\in\{1,\dots,\ell\}$ and assume that $y_{i_t}=m+pn_1+r$ for some $0\leq r<n_1$ and $0\leq p\leq n_0$. The case when $y_{i_t}=m+p(n_1-1)+(n_0+1)+r$ for some $0\leq r<n_1-1$ and $n_0<p< m$ is similar and will be omitted. Since $P\in \mbox{$_m\mathcal{D}_{m,n}(s)$}$ (hence $P$ is admissible) we must have
$$y_i-\ell_i>y_{i_t}-m+p,\ \ \text{$\forall i\in\{ i_t,\dots,i_{t+1}-1\}$}.$$
Equivalently,
\begin{equation*}y_{i_t}-p=m+p(n_1-1)+r>2y_{i_t}-y_i+\ell_i-m.\end{equation*}
Hence, if we write $2y_{i_t}-y_i+\ell_i-m=m+p'(n_1-1)+r'$ for some $0\leq r'<n_1-1$ and $0\leq p'\leq n_0$ we must have $p\geq p'$ which yields
\begin{equation}\label{admgenug}y_{i_t}-m>2y_{i_t}-y_i+\ell_i-2m+p\geq 2y_{i_t}-y_i+\ell_i-2m+p'.\end{equation}
We claim that \eqref{admgenug} already implies that $\Psi(P)$ is admissible for all pairs in $A(m,n-m)$ and hence $\Psi(P)\in \mathcal{D}_{m,n-m}(s)$. To see this, let $z\in \text{supp}(\Psi(P))$ such that $\ell_z(\Psi(P))\geq x_i+\ell_i-m$. We know that the path $\Psi(P)$ never goes below the line $y_{i_t}-m$ once it passes peak $B^{\Psi}_i$ and hence
\begin{equation}\label{admgenug2}\text{ht}_z(\Psi(P))\geq y_{i_t}-m.\end{equation}
Equation \eqref{admgenug} together with \eqref{admgenug2} proves that $\Psi$ is well--defined.
 It remains to verify that $\Psi$ is bijective. We define $\Psi^{-1}: \mathcal{D}_{m,n-m}(s-m) \rightarrow \mbox{$_m\mathcal{D}_{m,n}(s)$}$ as follows. Let $Q\in \mathcal{D}_{m,n-m}(s-m)$ a path with $d$ valleys whose coordinates are 
$$C_1=(x_1,y_1),\dots, C_d=(x_{d},y_{d}),\ \ x_1<x_2<\cdots< x_d.$$  
Furthermore, let $\zeta(Q)$ be the sequence from \eqref{incsubsetofvalleys}. We draw $v$ horizontal lines connecting the points $(x_{i_t},y_{i_t})$ with $(x_{i_t}-s_{i_t},y_{i_t})$, $1\leq t\leq v$, where $s_{i_t}=\min\{s\in \bn : (x_{i_t}-s,y_{i_t}) \in\text{supp}(Q)\}$ and reflect the path locally at these lines. Putting $m$ number of $1's$ in front of the reflected path yields $\Psi^{-1}(Q)$. The fact that $\Psi$ and $\Psi^{-1}$ are inverse to each other is clear and the fact that it is well--defined is similarly proven and we omit the details. 
\endproof
\end{prop}
\begin{cor}\label{eigenschaft2}
We have
$$_{m}\mathcal{E}_{s,n}^{1\rightarrow m}(q)=q^{(s-n)/2}\ _{1}\mathcal{E}_{s-m,n-m}^{1\rightarrow m}(q).$$
\proof
If $(i)\ (s-n)\notin 2\bz_+$ or $(ii)\ n<m$ or $(iii)\ m>s$, we have $_m\mathcal{D}_{m,n}(s)=\emptyset=\mathcal{D}_{m,n-m}(s-m)$ and the corollary is immediate. Otherwise the corollary follows from Proposition~\ref{eigenschaft1}.
\endproof
\end{cor}
\subsection{}

The following lemma will be needed.
\begin{lem}\label{neededlemma}Let $m,n,j,s\in \bn$ with $j<m$ and $s\geq j$. Then we have,
$$_{j}\mathcal{E}_{s,n}^{1\rightarrow m}(q)=\sum^{j}_{r=0} q^{(1-\delta_{r,0})(s-j)}\ \mbox{$_{(j+1-r-\delta_{r,s-j})}\mathcal{E}_{s-2r,n}^{1\rightarrow m}$}(q).$$
\proof
Each path $P\in \mbox{$_j\mathcal{D}_{m,n}(s)$}$ is of the form 
$$P=\underbrace{1\cdots 1}_{j}\underbrace{0\cdots 0}_{r}Q,$$
for some $r\in \{0,1\dots,j\}$ and a $01$--word $Q$, which is empty or starts with an up--step. Moreover, the path $Q$ is empty if and only if $r=s-j$. Define 
$$\tau(P):=\underbrace{1\cdots 1}_{j-r}Q\in\mbox{$_{(j+1-r-\delta_{r,s-j})}\mathcal{D}_{m,n}(s-2r)$}.$$ 
With other words we get a map
$$\tau:\mbox{$_j\mathcal{D}_{m,n}(s)$}\rightarrow \bigcup_{r=0}^{j} \mbox{$_{(j+1-r-\delta_{r,s-j})}\mathcal{D}_{m,n}(s-2r)$},$$
which is obviously well--defined and injective. 
\begin{figure}[H]
\begin{tikzpicture}[scale=0.3]
\draw [step=1,thin,gray!40] (0,0) grid (11,6);
\draw [very thick] (0,0) -- (3,3) -- (4,2)-- (5,3);
\draw [very thick,blue]  (5,3)--(6,4)--(8,2)--(11,5);
\node [left=1pt] at (1.5,5.5) {\scriptsize $P$};
\end{tikzpicture}
\hspace{1cm}
\begin{tikzpicture}[scale=0.3]
\draw [step=1,thin,gray!40] (0,0) grid (9,6);
\draw [very thick] (0,0) -- (3,3)--(4,4)--(6,2)--(9,5);
\node [left=1pt] at (2.8,5.3) {\scriptsize $\tau(P)$};
\end{tikzpicture}
\caption{$P=11101 {\color{blue}Q}$,\ \  $\tau(P)=111 Q$.}
\end{figure}
We show that $\tau$ is also surjective. If $r=s-j$ we have $|_{(2j-s)}\mathcal{D}_{m,n}(2j-s)|=\delta_{n,2j-s}$ and the preimage of the unique element in $_{(2j-s)}\mathcal{D}_{m,2js}(2j-s)$ is the path which consists of $j$ up--steps and $(s-j)$ down--steps. If $r\neq s-j$ choose $P'=\underbrace{1\cdots 1}_{j+1-r}Q'\in \mbox{$_{(j+1-r)}\mathcal{D}_{m,n}(s-2r)$}$ arbitrary, where $Q'$ is a $01$--word and assume without loss of generality that $s\geq 2r$ (c.f. \eqref{initial2}). Obviously the only candidate for the preimage is the path 
$$\tau^{-1}(P'):=\underbrace{1\cdots 1}_{j}\underbrace{0\cdots 0}_{r}Q,\ \ Q=1Q'.$$
We only have to show that $\tau^{-1}(P')$ is admissible. Note that the peaks of $P'$ and $\tau^{-1}(P)$ are identical if $r=0$ and else there is one additional peak at height $j$. Since $j<m$ there is no condition on the additional peak and $\tau^{-1}(P)$ is admissibe. Finally, it is clear that $\text{comaj}(P)=\text{comaj}(\tau(P))+(1-\delta_{r,0})(s-j)$, which finishes the proof.
\endproof
\end{lem}
\subsection{}\label{551}The proof of the main theorem requires one further representation theoretical result. For the rest of this section we consider only hook partitions $$\xi=(m,1^s),\ m\in\bn,\ s\in \bz_+$$ and recall from Remark~\ref{rem23} that 
\begin{equation}\label{fuseqloc}
V\left((1^s)\right)\cong W_{\text{loc}}(s).
\end{equation}
The proof of the next theorem is postponed to Section \ref{section6}.
\begin{thm}\label{sesnot} Let $s,m\in \bn$. Then, we have a short exact sequence of $\mathfrak{sl}_2[t]$--modules
$$0\rightarrow ker(\varphi)\rightarrow V\left((m,1^s)\right)\rightarrow V\left((m+1,1^{s-1})\right)\longrightarrow 0,$$
where  $ker(\varphi)$ can be filtered by 
$$\bigoplus_{r=\max\{m+2-s,1\}}^{m} \tau_s^* V\left((r,1^{s-2-m+r})\right) \oplus \tau_s^* M$$
and $M=\ev_0^* V_{\mathfrak{sl}_2}(m-s)$, if $s\leq m$ and $M=0$ otherwise.
\end{thm}
\begin{rem}
Another family of short exact sequences among fusion products has been constructed in \cite[Theorem 5(i)]{CV13} which coincides with the one given in Theorem~ \ref{sesnot} only if $m=1$. Generically they are different as all kernels given in \cite[Theorem 5(i)]{CV13} are proper fusion products.
\end{rem}
\subsection{Proof of Theorem~\ref{mainthm1}}From Proposition~\ref{existflag} we know that $V\left((m',1^s)\right)$ has a level $m$--Demazure flag if and only if $m\geq m'$. It is immediate from \eqref{fuseqloc} that the main theorem is a consequence of the following stronger statement, which gives a combinatorial model for graded multiplicities in fusion products associated to hook partitions.
\begin{thm}\label{hook}Let $m,m'\in\bn$, $s\in\bz_+$ such that $m\geq m'$. Let $\xi=(m',1^s)$. We have 
$$\mathcal{V}^{\xi\rightarrow m}_{n}(q)=\mbox{$_{m^{'}}\mathcal{E}_{s+m',n}^{1\rightarrow m}$}(q)=\sum_{P\in \hspace{0,03cm}  _{m'}\mathcal{D}_{m,n}(s+m')} q^{\text{comaj(P)}}.$$
\end{thm}
\proof
We prove the claim by induction on $s$. If $s=0$ we get 
$$\mathcal{V}^{\xi\rightarrow m}_{n}(q)=\delta_{m',n}=\mbox{$_{m^{'}}\mathcal{E}_{\xi,n}^{1\rightarrow m}$}(q)$$
and the induction begins. Assume that $s>0$. If $m>m'+1$ we can use Theorem \ref{sesnot} and our induction hypothesis to get
\begin{align}\mathcal{V}^{\xi\rightarrow m}_{n}(q)&=\notag \sum_{r=\max\{m'+2-s,1\}}^{m'+1}q^{(1-\delta_{r,m'+1})s}\hspace{0,03cm} \mathcal{V}^{(r,1^{s-2-m'+r})\rightarrow m}_{n}(q)+q^s\delta_{m,s+n}&\\&
=\label{needed2}\sum^{m'+1}_{r=1} q^{(1-\delta_{r,m'+1})s}\mbox{$_r\mathcal{E}_{s-2-m'+2r,n}^{1\rightarrow m}$}(q)+\mbox{$_{(m'-s)}\mathcal{E}_{m'-s,n}^{1\rightarrow m}$}(q).
\end{align}
Now Lemma \ref{neededlemma} implies that \eqref{needed2} is equal to \mbox{$_{m^{'}}\mathcal{E}_{s+m',n}^{1\rightarrow m}$} and the theorem is established in this case. So let $m=m'$ and note that Lemma~\ref{ausc} implies
$$\mathcal{V}^{\xi\rightarrow m}_{n}(q)=q^{(s+m-n)/2}\hspace{0,03cm}\mathcal{V}^{(1,1^{s-1})\rightarrow m}_{n-m}(q).$$
We apply the induction hypothesis to $\mathcal{V}^{(1,1^{s-1})\rightarrow m}_{n-m}(q)$ and conclude that
$$\mathcal{V}^{\xi\rightarrow m}_{n}(q)=q^{(s+m-n)/2}\mbox{$_1\mathcal{E}(s,m,n-m)$}.$$
The theorem is now immediate with Corollary \ref{eigenschaft2}.
\endproof

\section{Proof of Theorem~\ref{sesnot}}\label{section6}
\subsection{}
We first recall the more traditional definition of fusion products from \cite{FL99}. Let $n\in \bn$ and denote by $\ev_z^{*}V_{\mathfrak{sl}_2}(n)$, $z\in\bc$ the $\mathfrak{sl}_2[t]$--representation whose action is given by
$$(w\otimes f(t)).v=f(z) w.v,\ \ v\in V(n),\ f(t)\in \bc[t],\ x\in \mathfrak{sl}_2.$$
It is standard to show for pairwise distinct complex numbers $(z_1,\dots,z_{\ell})$ and a partition $\xi=(\xi_1\geq \cdots\geq \xi_{\ell})$ that $\ev_{z_1}^{*}V_{\mathfrak{sl}_2}(\xi_1)\otimes \cdots \otimes \ev_{z_{\ell}}^{*}V_{\mathfrak{sl}_2}(\xi_{\ell})$ is a cyclic representation for $\mathbf{U}(\mathfrak{sl}_2[t])$ with cyclic generator $v_{\xi}=v_{\xi_1}\otimes \cdots \otimes v_{\xi_{\ell}}$. Set 
$$\mathbf{U}(\mathfrak{sl}_2[t])_i=\text{span}\{(x_{i_1}\otimes f_1(t))\cdots(x_{i_{j}}\otimes f_{j}(t)): j\in \bz_+,\ \ \sum_{p=1}^{j}\text{deg }f_p(t)\leq i\}.$$
The associated graded space with respect to the filtration
$$0\subset \mathbf{U}(\mathfrak{sl}_2[t])_1v_{\xi}\subset \mathbf{U}(\mathfrak{sl}_2[t])_2v_{\xi}\subset \cdots $$
is isomorphic to $V(\xi)$ (see \cite[Theorem 5]{CV13}). This construction justifies why fusion products are sometimes called graded tensor products. Especially, we emphasize that the construction is independent of the chosen parameters and 
\begin{equation}\label{dimfor1}\text{dim } V(\xi)=\prod_{i=1}^{\ell} (\xi_i+1).\end{equation}
From the defining relations (see Section\ref{section23}) we immediately get the existence of a surjective homomorphism 
\begin{equation}\label{fusionsurjection}
\varphi: V\left((m,1^s)\right)\rightarrow V\left((m+1,1^{s-1})\right).
\end{equation} The proof of the next lemma is a straightforward calculation using \eqref{dimfor1}.
\begin{lem}\label{relaa}Let $\xi=(m,1^s)$. Then we have  $(y\otimes t^{s+1})v_{\xi}=0$ and 
$$\text{dim } ker \varphi=\sum _{r=\max\{m+2-s,1\}}^m \text{dim } V((r,1^{s-2-m+r}))+\max\{0,(m-s+1)\}.$$
\hfill\qed
\end{lem}
\subsection{}We need some more notation before we can prove Theorem~\ref{sesnot}. Let
\begin{equation*}_{\ell}\mathbf S(r,p)=\Big\{(b_k)_{k\ge 0}:b_k\in\mathbb Z_+,\ b_k=0\ \forall k<\ell \   \sum_{k\ge 0} b_k=r,\ \ \sum_{k\ge 0} kb_k=p\Big\},\ \ \mathbf S(r,p):=\mbox{$_0$}\mathbf S(r,p)\end{equation*}
and define
\begin{align*}_{\ell}\mathbf{y}(r,p)=&\sum_{\bob\in _{\ell}\mathbf S(r,p)}\ \prod_{i=0}^{p}(y\otimes t^i)^{(b_i)},\ \ \ \mathbf{y}(r,p):=\mbox{$_0$}\mathbf{y}(r,p).\end{align*}
We collect several immediate consequences.
\begin{lem}\label{simt1}
Let $V$ any representation of $\mathfrak{sl}_2[t]$ and $v\in V$, such that for all $a\in \bz_+$ we have $(h\otimes t^{a+1})v=(x\otimes t^a)v=0$. Then for all $r,p\in\bn$ we have 
\vspace{0,2cm}
\begin{enumerate}
\item $(x\otimes t)^{p}(y\otimes 1)^{r+p}v=\mathbf{y}(r,p)v$
\vspace{0,2cm}
\item 
$\mathbf{y}(r,p)-\mbox{$_{1}\mathbf{y}(r,p)$} \in \sum_{r'<r}\mathbf{U}(\mathfrak{sl}_2[t]) \mathbf{y}(r',p).$

\vspace{0,2cm}

\item Assume in addition $(y\otimes t^{p+1})v=0$. Then $(h\otimes t^k)\mathbf{y}(r,p)v$ is contained in the $\mathbf{U}(\mathfrak{sl}_2[t])$--span generated by
\begin{equation}\label{moduloset}
\{\mathbf{y}(r,p+k)v,\ \mathbf{y}(r-1,p')v: \  p'\geq p+1\}.\end{equation}
\end{enumerate}
\proof
The first claim of the lemma is a consequence of \cite[Lemma 7.1]{G78}. To prove part (2) we do induction on $r$. If $r=1$, we have
$$\mathbf{y}(1,p)-\mbox{$_{1}\mathbf{y}(1,p)$}=(y\otimes t^p)-(y\otimes t^p)=0.$$ If $r>0$, we use
$$\mathbf{y}(r,p)- \mbox{$_{1}\mathbf{y}(r,p)$} \in \sum_{r'<r}\mathbf{U}(\mathfrak{sl}_2[t])\mbox{$_{1}\mathbf{y}(r',p)$}$$
and apply our induction hypothesis to each $\mbox{$_{1}\mathbf{y}(r',p)$}$, which gives the desired result. The last claim is immediate if $k>p$; so assume that $k\leq p$. The claim  follows from the following calculation which we do modulo the $\mathbf{U}(\mathfrak{sl}_2[t])$--span of \eqref{moduloset}:
\begin{align*}(h\otimes t^k)\mathbf{y}(r,p)v&=\sum_{\bob\in \mathbf S(r,p+k)}(b_k+\cdots+b_p)\prod^{p}_{i=1}(y\otimes t^i)^{(b_i)}v&\\&=\sum_{\ell=0}^{k-1}\ \sum_{\bob\in \mathbf S(r,p+k)}(-b_{\ell})\prod^{p}_{i=1}(y\otimes t^i)^{(b_i)}v,\ \ \ \ \ \  (\text{subtract $r\mathbf{y}(r,p+k)v$})&\\&=\sum_{\ell=0}^{k-1}-(y\otimes t^{\ell})\sum_{\bob\in \mathbf S(r,p+k), b_{\ell}>0}(y\otimes 1)^{(b_0)}\cdots (y\otimes t^{\ell})^{(b_{\ell}-1)}\cdots (y\otimes t^{p})^{(b_{p})}v&\\&=\sum_{\ell=0}^{k-1}-(y\otimes t^{\ell})\ \mathbf{y}(r-1,p+k-\ell)v&\\&=0.\end{align*}
\endproof
\end{lem}
\subsection{}We are now able to determine generators of the kernel of the map \eqref{fusionsurjection}.
\begin{prop}\label{bestkern}Let $\xi=(m,1^s)$. The kernel of the homomorhism \eqref{fusionsurjection} is generated by the set 
\begin{equation}\label{genkernel}S=\{\mbox{$\mathbf{y}(r,s)$}v_{\xi}: 1\leq r\leq \min\{s,m\}\}.\end{equation}
\proof
Clearly Lemma \ref{simt1}(1) implies that $ker(\varphi)$ is generated by  
$$S':=\{\mbox{$\mathbf{y}(r,p)$}v_{\xi}: \ r,p\in \bn,\ r+p\geq 1+rk+\max\{0,(s-k)\}\ \ \text{for some $k\in \bn$}\}.$$ The claim follows if the prove that any element in $S'$ is contained in the $\mathbf{U}(\mathfrak{sl}_2[t])$--span of the elements in $S$. Let $\mathbf{y}(r,p)v_{\xi}\in S'$. If $r+p>1+rk+\max\{0,(s-k)\}$ we have 
$$r+p\geq 1+rk+\max\{0,(s-k)\}+1\geq  1+rk+\max\{0,(s+1-k)\}$$
and hence $\mbox{$\mathbf{y}(r,p)$}v_{\xi}=0$. So we can assume from now on that $r+p=1+rk+\max\{0,(s-k)\}$. If $k>1$ and $r\geq 2$ we get
$$r+p=1+r(k-1)+\max\{0,(s-k)\}+r\geq 1+r(k-1)+\max\{0,(s+1-(k-1))\},$$
and again $\mbox{$\mathbf{y}(r,p)$}v_{\xi}=0$. Otherwise if $k>1$ and $r=1$ we obtain
$$r+p=1+k+\max\{0,(s-k)\}\geq s+1\Rightarrow p\geq s.$$
This implies in the case when $p=s$ that $\mbox{$\mathbf{y}(1,p)$}v_{\xi}\in S$ and in the case when $p>s$ that $\mbox{$\mathbf{y}(1,p)$}v_{\xi}=0$ since
$$1+p\geq s+2\geq 1+k+\max\{0,(s+1-k)\}.$$ So if there exists an element $\mathbf{y}(r,p)v_{\xi}\in S'$ which is not contained in the $\mathbf{U}(\mathfrak{sl}_2[t])$--span of the elements in $S$ we must have $k=1$ and hence $p=s$. Moreover, if  $r>m$, then using Lemma \ref{simt1}(1) we get
$$r+s\geq m+s+1\Rightarrow (x\otimes t)^{(s)}(x\otimes t)^{(r+s)}v_{\xi}=\mbox{$\mathbf{y}(r,s)$}v_{\xi}=0.$$
So it remains to show that each element in 
$$\{\mbox{$\mathbf{y}(r,s)$}v_{\xi}: 1\leq r\leq m\}$$
is in the $\mathbf{U}(\mathfrak{sl}_2[t])$--span of the elements in $S$. If $s>m$ there is nothing to show; so assume from now on $s\leq m$.  If $r\leq s$, we already have that $\mathbf{y}(r,s)v_{\xi}\in S$, so we can assume in addition that $r>s$.
Let $(b_i)_{i\geq 0}$ a tuple of non--negative integers such that $$b_0+\cdots+b_s=r,\ \ b_1+2b_2+\cdots+sb_s=s.$$
We get \begin{align*}s&=b_1+\cdots+b_s+(b_2+2b_3+\cdots+(s-1)b_s)&\\&
=(r-b_0)+(b_2+2b_3+\cdots+(s-1)b_s)&\\&
\geq s+(r-s-b_0) \end{align*}
and hence $b_0\geq r-s>0$. This proves 
$$\mbox{$\mathbf{y}(r,s)$}v_{\xi}\in \sum_{b_0=r-s}^r \mathbf{U}(\mathfrak{sl}_2[t])\mbox{$_1\mathbf{y}(r-b_0,s)$}v_{\xi}.$$
The claim now follows from  Lemma\ref{simt1} (2).
\endproof
\end{prop}
\subsection{} An easy induction on $r$ shows that if $v$ is a weight vector of weight $N$ of a $\mathbf{U}(\mathfrak{sl}_2[t])$--representation $V$, then
$$[(x\otimes 1),(y\otimes 1)^{(r+1)}]v=(N-r)(y\otimes 1)^{(r)}v.$$
This implies together with Lemma \ref{moduloset}(1) that $(x\otimes 1)\mathbf{y}(r,s)v_{\xi}=A \mathbf{y}(r-1,s)v_{\xi}$ for some non--zero scalar $A\in\bc^*$. Hence using Proposition~\ref{bestkern} we get a filtration
$$0=V_0\subset V_1\subset V_2\subset \cdots \subset V_{\min\{s,m\}}=ker(\varphi),\ \ V_r:=\mathbf{U}(\mathfrak{sl}_2[t])\mathbf{y}(r,s)v_{\xi}.$$ 
Moreover,  we know from Lemma~\ref{relaa} that $(y\otimes t^{s+1})v_{\xi}=0$ and hence we obtain Lemma \ref{moduloset}(3)  that $(h\otimes t^k)\mathbf{y}(r,s)v_{\xi}$ is in the $\mathbf{U}(\mathfrak{sl}_2[t])$--span of the elements $\mathbf{y}(r,s+k)v_\xi,\mathbf{y}(r-1,s+k)v_\xi,\dots,\mathbf{y}(r-1,s+1)v_\xi$. Clearly the defining relations of $V(\xi)$ imply that all these elements are zero. Hence
$$(x\otimes 1)\mathbf{y}(r,s)v_{\xi}\in V_{r-1},\ \ (h\otimes t^k)\mathbf{y}(r,s)v_{\xi}\in V_{r-1},\ \ k\in \bn.$$
Since $(x\otimes \bc[t])\oplus (h\otimes t\bc[t])$ is generated by $\{(x\otimes 1), (h\otimes t^k): k\in \bn\}$ as a Lie algebra we obtain a well--defined map
\begin{equation}\label{brauchnur}W_{\text{loc}}(s+m-2r)\rightarrow V_r/V_{r-1},\ 1\leq r \leq \min\{s,m\}.\end{equation}
The following proposition along with Lemmn\ref{relaa} finishes the proof of Theorem~\ref{sesnot}.
\begin{prop}For $1\leq r< \min\{s,m\}$ there exists a surjective homomorphism
$$V\left((m+1-r,1^{s-1-r})\right)\rightarrow V_r/V_{r-1}.$$
Moreover, we have a surjective map $W_{\text{loc}}(s-m)\rightarrow V_{m}/V_{m-1}$ if $s\geq m$ and otherwise $$V_{s}/V_{s-1}\cong 
\ev_0^* V_{\mathfrak{sl}_2}(m-s).$$  
\proof Assume that $1\leq r < \min\{s,m\}$. From \eqref{brauchnur} and the defining relations of fusion products it will be enough to show that the following holds in $V_r/V_{r-1}$:
$$(x\otimes t)^{(p')}(y\otimes 1)^{(r'+p')}\mbox{$_1\mathbf{y}(r,s)$}v_{\xi}=0,\ \forall r',p',k'\in \bn: r'+p'\geq 1+r'k'+\max\{0,(s-r-k')\}.$$ 
Recall the traditional construction of fusion products from Section~\ref{551} and the independence of the choice of the parameters. Without loss of generality we will identify $v_{\xi}$ with $v_m\otimes v_1\otimes \cdots \otimes v_1$ and choose the first evaluation parameter to be zero, i.e. $z_1=0$. We get
\begin{align*}(x\otimes t)^{(p')}(y\otimes 1)^{(r'+p')}\mbox{$_1\mathbf{y}(r,s)$}v_{\xi}&=(x\otimes t)^{(p')}(y\otimes 1)^{(r'+p')}\left(v_{m}\otimes \mbox{$_1\mathbf{y}(r,s)$}(v_1\otimes \cdots \otimes v_1)\right)&\\&=v_{m}\otimes \left((x\otimes t)^{(p')}(y\otimes 1)^{(r'+p')}\mbox{$_1\mathbf{y}(r,s)$}(v_1\otimes \cdots \otimes v_1)\right).\end{align*}
Since $r'+p'+r>s$ and the lowest weight of $V_{\mathfrak{sl}_2}(1)^{\otimes s}$ is $(-s)$ we must have 
$$ \left((x\otimes t)^{(p')}(y\otimes 1)^{(r'+p')}\mbox{$_1\mathbf{y}(r,s)$}(v_1\otimes \cdots \otimes v_1)\right)=0.$$
It remains to consider the case $r=\min\{s,m\}$. If $m\geq s$ we obtain with Lemma \ref{moduloset}(2) that the following holds in $V_s/V_{s-1}$:
$$(y\otimes t^k)\mathbf{y}(s,s)v_{\xi}=(y\otimes t^k)\mbox{$_1\mathbf{y}(s,s)$}v_{\xi}=(y\otimes t^k)(y\otimes t)^{(s)}v_{\xi}=0,\ \forall k\in \bn.$$
Hence $V_s/V_{s-1}\cong \ev_0^* V_{\mathfrak{sl}_2}(m-s)$. If $s\geq m$, the claim follows already from the discussion preceding the proposition, see \eqref{brauchnur}.
\endproof
\end{prop}
\bibliographystyle{plain}
\bibliography{bibfile}
\end{document}